\documentclass[11pt,a4paper, twoside]{article}
\usepackage[notref,notcite]{showkeys}

\usepackage{amsmath,amssymb,amsthm,amsfonts,mathrsfs,amscd,environ}
\usepackage{dsfont}
\usepackage{latexsym,enumerate,color,geometry,extarrows}
\usepackage{xcolor}
\usepackage{verbatim,fancyhdr,enumitem}
 \NewEnviron{ews}{%
\begin{equation}\begin{split}
  \BODY
\end{split}\end{equation}
}

\NewEnviron{ews*}{%
\begin{equation*}\begin{split}
  \BODY
\end{split}\end{equation*}
}

\def\beg{\begin}
\def\bequ{\begin{equation}}
\def\enqu{\end{equation}}
\def\bes{\begin{split}}
\def\ens{\end{split}}
\def\bews{\begin{ews}}
\def\beqn{\begin{eqnarray}}
\def\enqn{\end{eqnarray}}
\def\beq*{\begin{equation*}}
\def\enq*{\end{equation*}}
\def\bqn*{\begin{eqnarray*}}
\def\eqn*{\end{eqnarray*}}
\def\bary{\begin{array}}
\def\eary{\end{array}}
\def\bpma{\begin{pmatrix}}
\def\epma{\end{pmatrix}}
\def\bvma{\begin{Vmatrix}}
\def\evma{\end{Vmatrix}}

 \numberwithin{equation}{section}

\def\al{\alpha}
\def\be{\beta}
\def\ga{\gamma}
\def\de{\delta}
\def\ep{\epsilon}

\def\et{\eta}
\def\th{\theta}

\def\ka{\kappa}
\def\la{\lambda}

\def\si{\sigma}
\def\ta{\tau}

\def\De{\Delta}

\def\Ph{\Phi}

\def\Om{\Omega}

\def\R{\mathbb R}
\def\P{\mathbb P}
\def\E{\mathbb E}

\def\N{\mathbb N}

\def\sE{\mathscr E}
\def\sF{\mathscr F}
\def\sD{\mathscr D}
\def\sC{\mathscr C}
\def\sB{\mathscr B}

\def\sL{\mathscr L}

\def\sP{\mathscr P}

\def\cT{\mathcal T}
\def\cM{\mathcal M}

\def\d{\mathrm{d}}

\def\ff{\frac}
\def\ra{\rightarrow}
\def\nn{\nabla}
\def\pp{\partial}
\def\<{\langle}
\def\>{\rangle}
\def\sq{\sqrt}
\def\tld{\tilde}
\def\we{\wedge}
\def\1{\mathds{1}}

\def\trac{\mathrm{tr}}

\def\Ent{\mathrm{Ent}}

\allowdisplaybreaks

\title{{\bf Existence and non-uniqueness of  stationary distributions  for distribution dependent SDEs}
}

\author{
{\bf Shao-Qin Zhang }\\
\footnotesize{School of Statistics and Mathematics, Central University of Finance and Economics, Beijing 100081, China}\\
\footnotesize{Email: zhangsq@cufe.edu.cn}\\
}

\begin{document}

\maketitle

\begin{abstract}
The existence of stationary distributions to distribution dependent stochastic differential equations are investigated by using the ergodicity of the associated decoupled equation and the Schauder fixed point theorem. By using Zvonkin's transformation, we also establish  the existence result for equations with  singular coefficients.  Instead of the uniqueness, the non-uniqueness of stationary distributions are considered for equations with regular coefficients. Concrete examples  including McKean-Vlasov stochastic equations with the quadratic interaction and the  non-quadratic  interaction, and equations with a bounded and discontinuous drift are presented to illustrate our non-uniqueness results. 
\end{abstract}\noindent

\noindent


\vskip 2cm

\section{Introduction}

When investigating the propagation of chaos  for  interacting diffusions, McKean \cite{McK} introduced a nonlinear stochastic differential equation (SDE) whose coefficients depend  on the own law of the solution. This SDE is referred  as the McKean(-Vlasov) SDE, established from systems of interacting diffusions by passing to the mean field limit. The associated empirical measures  converge in the weak sense  to a probability measure with density which satisfies a nonlinear parabolic partial differential equation (PDE) called McKean-Vlasov equation in the literature, see e.g. \cite{Mel,Szn}.  The existence of several stationary distributions to the McKean-Vlasov equation can occur for non-convex the confining potential.   Many works are concerned about this phenomenon (referred as phase transition)  for McKean-Vlasov equations. For instance, \cite{Tam} provided a  criteria for McKean-Vlasov equations with  odd interaction potentials; \cite{Daw} established the phase transition for the equation  with a particular double-well confinement and Curie-Weiss interaction on the line; quantitative analyses for continuous and discontinuous phase transition were provided for McKean-Vlasov equations on the torus  in \cite{CGPS,ChPa}; equations with multi-wells landscape on the whole space were investigated extensively  by Tugaut et  al. in \cite{DuTu,HT10a,HT10b,Tug10,Tug13,Tug14a,Tug14b}. There exist analogues called nonlinear master equations to describe nonlinear pure jump Markov processes. The phase transition is also studied for  nonlinear master equations, see e.g. \cite{Chen,FenZ}.  In this paper, we aim to investigate the existence and non-uniqueness of stationary distribution   to distribution dependent SDEs which cover McKean-Vlasov SDEs as examples.

Let $\sP $ be the space of probability measures on $\R^d$ equipped  with the weak  topology,  $W$ be a $d$-dimensional Brownian motion on a complete filtration probability space $(\Om,\sF,\{\sF_t\}_{t\ge 0},\P)$, and let $\sL_\et$ be the law of the random variable $\et$. We consider the  SDE on $\R^d$ of the following form
\beg{align}\label{equ00}
\d X_t=b(X_t,\sL_{X_t})\d t+\si(X_t,\sL_{X_t})\d W_t,
\end{align}
where the coefficients 
$$b:~\R^d\times \sP \ra \R^d,\qquad \si:~\R^d\times\sP \ra \R^d\otimes\R^d$$
are measurable. Solutions to \eqref{equ00} in strong and weak sense are extensively studied by many works, see e.g. \cite{HW,HRW,Wan18} and references therein. In this paper, we aim to find  stationary distribution   to \eqref{equ00}, i.e. $\mu\in\sP $ so that for $\sL_{X_0}=\mu$, there is a solution   with $\sL_{X_t}\equiv \mu$ for all $t\geq 0$.  The solution $X_t$ associated to the stationary probability measure $\mu$ is called the stationary solution. For $b,\si$ are independent of the distribution of the solution, $\mu$ is the invariant probability of the associated linear Markov process.  

By setting $b(x,\mu)=-\nn V(x)-\nn F*\mu(x)$, we get the McKean-Vlasov SDE with that the confining potential $V$ and  the interaction potential $F$ are  differentiable functions on $\R^d$, where $*$ stands for the convolution on $\R^d$:
$$f*\mu(x)=\int_{\R^d} f(x-y)\mu(\d y),~f\in\sB(\R^d).$$
Moreover, in the case that $\si=\si_0 I$ for some $\si_0\in\R$,  stationary distributions  are of an explicit formulation:
\beg{align}\label{eq-st-exp}
\mu(\d x)=\ff {\exp\left\{-\ff  2 {\si_0^2} \left( V(x)+ F*\mu(x) \right)\right\}} {\int_{\R^d} \exp\left\{-\ff  2{\si_0^2} \left( V(x)+ F*\mu(x) \right)\right\} \d x}\d x.
\end{align}
Then  the  existence of stationary measures can be investigated by this explicit formulation and the fixed point theorem for  multi-well confinement $V$  and  even polynomial interaction  $F$  (see e.g. \cite{Daw,DuTu,Tug14a}), or combining it with the free energy functional associated to the McKean-Vlasov equation (see \cite{CGPS,ChPa,Tug13,Tug14b} for instance): 
\beg{align*}
\sE^{V,F}(\mu):=\ff {\si_0^2} 2 Ent(\mu|\mu_{V_{\si_0}}) +\ff 1 2\mu(F*\mu),
\end{align*}
where $\mu_{V_{\si_0}}(\d x)=Z_0 \exp\{- 2 \si_0^{-2}V(x)\}\d x$ is a probability measure with  the normalizing constant $Z_0$ and $\Ent(\mu|\mu_{V_{\si_0}})$ is the classical relative entropy. However, for SDE \eqref{equ00}, the explicit formulation as \eqref{eq-st-exp}  for stationary measures is not available,  and the free energy functional $\sE^{V,F}$ may be not available or  less explicit even it exists.  

For distribution dependent SDEs, the existence and uniqueness of the invariant probability measures  and the exponential convergence to the unique invariant probability measure  in the Wasserstein (quasi-)distance or the entropy have been investigated, see e.g. \cite{EGZ,Mal,RenWa,Wan18,Wan21} and so on.  Results on the exponential convergence  can lead to the existence of a contractive mapping, see e.g. \cite{Wan18}. However,  the contracting-mapping principle implies the existence and uniqueness at once time.   This excludes stochastic equations  with  several invariant measures.  In this paper, we use the ergodicity of a decoupled SDE (\eqref{equ0} below) and the Schauder fixed point theorem  instead of the  contractive-mapping principle to find  stationary distributions  of \eqref{equ00}. Given $\mu\in\sP(\R^d)$, we get from \eqref{equ00}  the following decoupled SDE by freezing $\sL_{X_t}\equiv \mu$
\beg{align}\label{equ0}
\d X_t^\mu=b(X_t^\mu,\mu)\d t+\si(X_t^\mu,\mu)\d W_t.
\end{align}
Define $P_t^\mu f(x)=\E f(X_t^{\mu,x})$ the associated Markov semigroup to \eqref{equ0}. If $P_t^\mu$ has a unique  invariant probability measure, then   a mapping $\cT$ on $\sP $ with $\cT_\mu$ be the invariant probability measure of $P_t^\mu$ is established. Thus the existence of stationary measures to \eqref{equ00} can be reformulated as the existence of  fixed points of $\cT$.  The existence results are given for SDEs with regular coefficients and with singular  coefficients, see Theorem \ref{thm1} and Theorem \ref{thm2} below. Non-uniqueness of stationary distributions  are investigated and some criterion are given, see Theorem \ref{thmN} in Section 2.

This paper is structured as follows. Section 2 is devoted to the existence and non-uniqueness results on  stationary distributions of SDEs with regular coefficients,  and concrete examples are given in this section to illustrate the non-uniqueness of stationary distributions. We extends the existence result to  stochastic equations with singular coefficients in Section 3, and an example for SDEs with a bounded measurable drift is presented.

\section{Existence and non-uniqueness results for  regular SDEs}
In this section, we consider \eqref{equ00} with the coefficients $b$ and $\si$  are regular.  Let 
\beg{align*}
\sP^r&=\{\mu\in\sP ~|~\|\mu\|_r:=(\mu(|\cdot|^r))^{\ff 1 r}<\infty\}\\
\sP^r_M&=\{\mu\in\sP ~|~\|\mu\|_r:=(\mu(|\cdot|^r))^{\ff 1 r}\leq M\},~r>0,M>0.
\end{align*}
Denote by $\|\cdot\|$ the operator norm of a matrix and $\|\cdot\|_{HS}$ the Hilbert-Schmidt norm. For any measurable  matrix-value function $f$ on $\R^d\times\sP$, denote 
$$ \|f\|_\infty= \sup_{x\in\R^d,\mu\in\sP}\|f(x,\mu)\|_{HS}.$$

Assume that $b$ and $\si$ satisfy the following hypothesis:
\beg{description}[align=left, noitemsep]
\item[(H1)] There exist constants $r_1\geq 0,r_2\geq r_1,r_3>0$,  $C_1 >0$, and nonnegative   $C_2,C_3$ such that for any $\mu\in\sP^{1+r_2} $  
\beg{align}\label{ine-b-si0}
 2\<b (x,\mu),x\>+\|\si(x,\mu)\|_{HS}^2&\leq -C_1|x|^{1+r_1}+C_2+C_3\|\mu\|_{1+r_2}^{r_3}.
\end{align}
\item[(H2)]  For every $n\in\N$ and $\mu\in\sP^{1+r_2}$, there exists  $K_n>0$ such that 
\beg{align}\label{b0-K}
|b(x,\mu)-b(y,\mu)| &+\|\si(x,\mu)-\si(y,\mu)\|_{HS}\leq K_n|x-y|,~|x|\vee |y|\leq n.
\end{align}
There exist non-negative constants $C_4$, $C_5$ and $C_6$ such that
\beg{align}\label{ine-b-si00}
\|\si(x,\mu)\|_{HS}^2\leq C_4|x|^{1\we r_1+r_1}&+C_5\|\mu\|_{1+r_2}^{r_3}+C_6,~x\in\R^d,\mu\in\sP^{1+r_2}.
\end{align}
For every $\mu\in\sP^{1+r_2}$, there exists $C_7>0$ such that
\beg{align}\label{ineb-g}
|b(x,\mu)|\leq C_7(1+|x|^{r_1}),~x\in\R^d.
\end{align}
\item[(H3)]  For each $n\geq 1$, and $\mu_m,\mu\in\sP^{1+r_2}_M$ with $\mu_m\xrightarrow{w}\mu$, there is 
\beg{align}\label{ineb}
\lim_{  m\ra+\infty }\sup_{|x|\leq n}\left(|b(x,\mu)-b(x,\mu_m)|+\|\si(x,\mu)-\si(x,\mu_m)\|_{HS}\right)=0.
\end{align}
\end{description}
Then we have the following result on the existence of  stationary distributions  for \eqref{equ00}. In \eqref{non-de} below or more generally for a symmetric matrix  $A\in\R^d\otimes\R^d$ and $\la\in\R$, the inequality $A\geq \la$  means that
\beg{align*}
\<Av,v\>\geq \la|v|^2,~v\in\R^d,
\end{align*}
and $A>\la$ is defined similarly. 

\beg{thm}\label{thm1}
Assume (H1)-(H3) and that $\si$ is non-degenerate:
\beg{align}\label{non-de}
\si(x,\mu)\si^*(x,\mu)>0,~x\in\R^d,\mu\in\sP^{1+r_2}.
\end{align}
If $C_1$, $C_3$, $C_4$, $C_5$, and $r_1,r_2,r_3$ satisfy
\beg{align}\label{ine-C-r}
&r_3 \leq   1+r_1  ;\qquad  C_1-(r_2-r_1)C_4>0;\\
&\ff {C_3+C_5(r_2-r_1)} {C_1-C_4(r_2-r_1)}<1,~\text{when}~r_3 = 1+r_1,\label{ine-C-r'}
\end{align}
then \eqref{equ00} has a stationary distribution. 
\end{thm}

Set $b(x,\mu)=-\nn V(x)- \nn F*\mu$ for some twice continuous differential functions $V$ and $F$. Then we have the following corollary which will be used in discussing the non-uniqueness of stationary distributions .  
\beg{cor}\label{cor0}
Assume that $\si$ is non-degenerate, bounded on $\R^d\times\sP$ and satisfies (H2) and (H3), and $V$, $F$ are twice continuous differential functions with  non-negative constants  $\al_0,\al_1,\al_2,\al_3$ and $\ga_0\in [0,3)$ such that
\beg{align}\label{nnV1}
|\nn V(x)|&\leq \al_0(1+|x|^3),\\
|\nn F(x)|&\leq \al_1+\al_2|x|^3,\label{nnF1}\\
\|\nn^2 F(x)\|&\leq \al_3(1+ |x|^{\ga_0}),~x\in\R^d. \label{n2F4}
\end{align}
Suppose that there exist constants $\be_0,\be_1$, and $\be_2>0$, $\be_3\geq 0$ such that 
\beg{align}\label{VF0}
\nn^2 V(x)+\nn^2 F(x-y)&\geq \be_0-2\be_1|x|+3\be_2|x|^2-\be_3|y|^2,~x,y\in\R^d,
\end{align}
and  $\al_2+\be_3<\be_2$. Then \eqref{equ00} has stationary distributions  in $\sP^4$.

\end{cor}

Let $a\in\R^d$. Denote by $\mu_a$ the shifted probability  of $\mu$ by $a$: 
$$\mu_a(f)=\mu(f(\cdot-a)).$$ 
The next theorem is devoted to a suffice condition to find a stationary probability  measure in  sets of following form where $\ka_2\geq \ka_1>0$, and $0<\ga_1\leq \ga_2<1+r_2$, 
$$\sP_{a, \ka_1,\ka_2}^{ \ga_1,\ga_2}=\{\mu\in\sP^{1+r_2}  ~|~\|\mu_a\|_{\ga_1} \leq \ka_1 ,\|\mu_a\|_{\ga_2}\leq \ka_2 \}.$$
Then  non-uniqueness can be established by using this  criteria.  Denote
$$\sD_{\ka_1,\ka_2}=\{(r_1,r_2)~|~0\leq r_1\leq r_2\leq \ka_2,0\leq r_1\leq \ka_1\}.$$ 

\beg{thm}\label{thmN}
Assumptions of Theorem \ref{thm1} hold. Assume that there are $a\in\R^d$, $\ga_1,\ga_2\in (0,1+r_2)$ with $\ga_1\leq \ga_2$, and measurable functions $g_1,g_2$ on $[0,+\infty)^3$ so that $g_1(\cdot,r_1,r_2)$ and  $g_2(\cdot,r_1,r_2)$ are continuous and convex on $[0,+\infty)$ and  for  $\mu\in\sP^{1+r_2}$ 
\beg{align}\label{bgg}
&2\<b(x+a,\mu),x\>+\|\si(x+a,\mu)\|_{HS}^2\nonumber\\
&\qquad\qquad \leq -g_1(|x|^{\ga_1},\|\mu_a\|_{\ga_1},\|\mu_a\|_{\ga_2})  -g_2(|x|^{\ga_2},\|\mu_a\|_{\ga_1},\|\mu_a\|_{\ga_2}),~x\in\R^d.
\end{align}
Let   
$$g(r,r_1,r_2) =g_1(r,r_1,r_2)+g_2(r^{ \ff {\ga_2} {\ga_1} },r_1,r_2),~r_1,r_2\geq 0.$$
Suppose  that $g(\cdot,r_1,r_2)$ is convex on $[0,+\infty)$ and  
there exist  $\ka_2\geq \ka_1>0$  such that for all $(r_1,r_2)\in\sD_{\ka_1,\ka_2}$, 
\beg{align}\label{sset1}
&g(r^{\ga_1},r_1,r_2) > 0, \quad r\in [\ka_1,+\infty),\\
g_2(r^{\ga_2},r_1,r_2) & + \inf_{r\in [0,\ka_1]} g_1(r^{\ga_1},r_1,r_2)>0,\quad r\in [\ka_2,+\infty).\label{sset2}
\end{align} 
Then there is $\mu\in\sP_{a, \ka_1,\ka_2 }^{\ga_1,\ga_2}$ being a stationary probability measure of \eqref{equ00}.  

Consequently, if there exist  $a_1,a_2\in\R^d$ and $\ka_1 <  \ff {|a_1-a_2|} 4 $ such that the above assumptions  hold, then \eqref{equ00} has two distinct  stationary probabilities $\mu_1\in \sP_{a_1, \ka_1,\ka_2 }^{\ga_1,\ga_2},\mu_2\in \sP_{a_2, \ka_1,\ka_2 }^{\ga_1,\ga_2}$. 
\end{thm}

If $\ga_1=\ga_2=:\ga$, we denote $\sP_{a,\ka}^{\ga}=\sP_{a, \ka ,\ka }^{ \ga ,\ga }$ and have simplified criteria as follows. 
\beg{cor}\label{cor1}
Assumptions of Theorem \ref{thm1} hold. Assume that there are $a\in\R^d$, $\ga \in (0,1+r_2)$ and  a measurable function $g$ on $[0,+\infty)^2$ so that $g(\cdot,r_1 )$ is continuous and convex  for each $r_1\geq 0$,  and  
\beg{align}\label{bgg1}
2\<b(x+a,\mu),x\>+\|\si(x+a,\mu)\|_{HS}^2\leq -g(|x|^{\ga },\|\mu_a\|_{\ga }).
\end{align}
If there exist  $\ka>0$  such that 
\beg{align}\label{sset3}
g (r^{\ga },r_1 )>0,~r\geq \ka, 0\leq r_1\leq \ka,
\end{align}
then there is $\mu\in\sP_{a, \ka}^{ \ga}$ being a stationary probability measure of \eqref{equ00}.  Consequently, if there exist  $a_1,a_2\in\R^d$ and $\ka <  \ff {|a_1-a_2|} 4 $ such that the above assumptions  hold, then \eqref{equ00} has two different stationary probabilities $\mu_1\in \sP_{a_1, \ka }^{ \ga},\mu_2\in \sP_{a_2, \ka}^{\ga}$.
\end{cor}

For the McKean-Vlasov SDEs, we have the following corollary.  

\beg{cor}\label{cor2}
Set $b(x,\mu)=-\nn V(x)-\nn F*\mu(x)$ for twice continuous differentiable functions $V$ and $F$.  Assume that $\nn V$ has polynomial growth, and that  $\si$ is non-degenerate, bounded on $\R^d\times\sP$ and satisfies (H2) and (H3). The point  $a\in\R^d$ is a critical point of   $V$, and there are positive constants $\be_0,\be_2$  and non-negative constants $\be_1$  such that  for every $x,y\in\R^d$ 
\beg{align}\label{VF}
 \nn^2 V(a+x)+\nn^2 F(a-y+x) \geq  \be_0-2\be_1|x|+3\be_2|x|^2.
\end{align}
Suppose there are non-negative constants $\al_0,\al_1,\al_2,\al_3 $ such that  \eqref{nnV1} and \eqref{n2F4} with $\ga_0=1$ hold, and 
\beg{align}\label{LipF}
|\nn F(x)|&\leq \al_1+\al_2|x|,~x\in\R^d.
\end{align}
Let
\beg{align*}
\ka_0 &= \inf\{r>0~|~\be_2r^3-\be_1r^2+(\be_0-\al_2)r-\al_1>0\},\\
\ka_1 &= \inf\{r>\ka_0~|~\be_2r^3-\be_1r^2+(\be_0-\al_2)r-\al_1<0\},\\
\ka_2 &= \sup\{r>\ka_0~|~\be_2r^3-\be_1r^2+(\be_0-\al_2)r-\al_1<0\}.
\end{align*}
If $\be_0>\al_2$, $\be_0\be_2>3\be_1^2/8$ and there is $\ka\in (\ka_0,\ka_1)\cup (\ka_2,+\infty)$  such that
\beg{align}\label{si-be-al}
\|\si\|_\infty^2< 2\left(\be_2 \ka^3-\be_1 \ka^2+ (\be_0-\al_2) \ka - \al_1  \right)\ka,
\end{align}
then  there is a stationary probability measure $\mu\in\sP^{1}_{a, \ka}$.
\end{cor}
\beg{rem}
The inequality \eqref{VF} may hold if $a\in\R^d$ is a minimum of $V$.  In the following examples,  stationary distributions  are found around the minimums of $V$, see also \cite{Tug14a}.

Corollary \ref{cor2} provides a sufficient condition to find the stationary probability measure around $a$. According to Theorem \ref{thmN}, to prove that there are several stationary distributions, $\ka$ needs to be small (less than a quarter of the distance between two minimums). To get small $\ka$ such that \eqref{si-be-al} holds, a sufficient condition is that $\al_1=0$ and $\|\si\|_\infty$ is small, see examples below. 
\end{rem}
 
We present concrete  examples to illustrate the non-uniqueness of the stationary distributions.
\beg{exa}\label{ex-0}
Let $d=1$, $a_1,a_2\in\R$ with $a_1a_2<0$,  $\be>0$ and $\al>0$.   Consider the following McKean-Vlasov SDE with quadratic interaction
\beg{align}\label{eq-ex1}
\d X_t&=-\be(X_t-a_1)X_t(X_t-a_2)\d t-\al \int_{\R}(X_t-y)\sL_{X_t}(\d y)\d t\nonumber\\
&\quad +\si(X_t,\sL_{X_t})\d W_t.
\end{align}
Assume that $\si$ is positive and bounded on $\R\times\sP$  and satisfies \eqref{b0-K} and (H3).  Then \eqref{eq-ex1} has a stationary distribution. Moreover, 
For $a_1,a_2,\al,\be$ satisfying 
\beg{align}\label{aaba}
a_1^2+a_2^2+2(a_1-a_2)^2+a_1^2\vee a_2^2<\ff {8\al} {\be}
\end{align}
and for some $\ka\in (0, (|a_1|\we|a_2|)/2)$, 
\beg{align}\label{in-sika}
\|\si\|_\infty< \ka\sq{2\be(\ka-|a_1-a_2|)(\ka-|a_1|\we |a_2|)},
\end{align} 
there exist two distinct stationary distributions  $\nu_1,\nu_2\in\sP^{1+r_2}$ such that
\beg{align}\label{nunu}
\nu_1(|\cdot-a_1|)\leq \ka,\qquad \nu_2(|\cdot-a_2|)\leq \ka.
\end{align}
Consequently, if $\si$ is a positive constant,  $a_1=-a_2$ with \eqref{aaba} holds, and $\|\si\|_\infty<\ka\sq{2\be(\ka-2|a_1|)(\ka-|a_1|)}$ for some $0<\ka<|a_1|/2$, then \eqref{eq-ex1} has at least three stationary distributions.
\end{exa}

The following example shows that our criteria can also be used to deal with McKean-Vlasov SDEs with non-quadratic interaction.

\beg{exa}\label{ex-1}
Let $a_1,a_2\in\R^d$, and let $\be $, $\al_1$ and $\al_2$ be positive constants. Consider the following McKean-Vlasov SDE with non-quadratic interaction
\beg{align}\label{eq-ex2}
\d X_t&=-\ff {\be} 2\left\{(X_t-a_1) |X_t-a_2|^2+(X_t-a_2)|X_t-a_1|^2\right\}\d t+\si(X_t,\sL_{X_t})\d W_t\nonumber\\
&\qquad - \int_{\R}(\al_1|X_t-y|^2(X_t-y)+\al_2 (X_t-y))\sL_{X_t}(\d y)\d t.
\end{align}
Assume that  $\si$ is positive and bounded on $\R\times\sP$  and satisfies \eqref{b0-K} and (H3). Let 
\beg{align*}
\th_0 &=\ff {\al_2} {\be |a_1-a_2|^2},\qquad \th_1=\ff {(\be+ {\al_1} )^{\ff 3 2} \be^{\ff 1 2} } {2(3(4+\th_0))^{\ff 3 2}\al_1^2},\\
\ka_1&=\th_1 |a_1-a_2|,\qquad \ka_2=\left(\ff {\be\ka_1} {4\al_1}\right)^{\ff 1 3} |a_1-a_2|^{\ff 2 3}.
\end{align*}
Suppose that
\beg{align}  \label{aba1a2} 
 \al_1  + {\be}   & >  243(4+\th_0) {\be}, \\
\ff {\|\si\|_\infty^2} { 2\be|a_1-a_2|^4} & < \left( \ff {\th_1(4+\th_0)} 2   -\ff {3\al_1} {\be}\right)\we \ff { \th_1^2} 5 . \label{si-ka12}
\end{align}
Then there exist two distinct stationary distributions  $\nu_1,\nu_2\in\sP^{1+r_2}$ such that
\beg{align}\label{nunu}
\nu_i(|\cdot-a_i|)\leq \ka_1,\qquad \left(\nu_i(|\cdot-a_2|^3)\right)^{\ff 1 3}\leq \ka_2,~i=1,2.
\end{align}
\end{exa}

The proof of Theorem \ref{thm1} is given in the following subsection, and  proofs of Theorem \ref{thmN}  is shown in Subsection 2.2. This section will ends up with the proofs of examples in Subsection 2.3.

\subsection{Proof of Theorem \ref{thm1}}
To prove Theorem \ref{thm1}, we prove firstly that $\cT$ is well-defined. The space $\sP^1$ equipped with the Kantorovich-Rubinstein-Wasserstein distance ($W$-distance for short) is a complete metric space (see e.g. \cite[Theorem 5.4]{Chen}):
\beg{align*}
W(\mu,\nu):= \inf_{\pi\in\sC(\mu,\nu)}\int_{\R^d\times\R^d}|x-y|\pi(\d x,\d y),~\mu,\nu\in\sP^1,
\end{align*} 
where $\sC(\mu,\nu)$ consists of all the couplings of $\mu$ and $\nu$. Let $\cM^1$ be the set of all finite signed measures on $\R^d$ with $|\mu|(|\cdot|)<\infty$, and let 
$$\|\mu\|_{KR}:=|\mu(\R^d)|+\sup_{h\in Lip(\R^d),h(0)=0} \int_{\R^d} h(x)\mu(\d x),~\mu\in\cM^1.$$
Then $(\cM^1,\|\cdot\|_{KR})$ is a normed space. Moreover,   
\beg{align*}
\|\mu-\nu\|_{KR}=W(\mu,\nu),~\mu,\nu\in\sP^1,
\end{align*}
see \cite[Corollary 5.4]{CarDe}. This, together with that $(\sP^1,W)$ is a complete metric space,  yields that the Schauder fixed point theorem (see e.g. \cite[Theorem 8.8]{Deim} is available on a nonempty, convex, and compact  subset of $(\sP^1,W)$. Hence,  we also need to prove that $\cT$ is continuous in some suitable compact subset of $(\sP^1,W)$.

The  lemma below gives a fundamental estimate on $X_t^{\mu,x}$ which will be used   later.
\beg{lem}\label{lem1}
Assume that (H1) holds and 
\beg{align}\label{inesi}
\|\si(x,\mu)\|_{HS}^2\leq C_4|x|^{1+r_1}+C_5\|\mu\|_{1+r_2}^{r_3}+C_6,~x\in\R^d, \mu\in\sP^{1+r_2}.
\end{align}
Then 
\beg{align}\label{ine-lem2}
\E\sup_{t\in [0,T]}|X_t^{\mu,x}|^2&\leq \left(2+\ff {  C_4} {C_1}\right) \left(|x|^2+\left(C_2+C_3\|\mu\|_{1+r_2}^{r_3}\right)T\right)\nonumber\\
&\qquad + \left(C_5\|\mu\|_{1+r_2}^{r_3}+C_6\right)T,~T\geq 0.
\end{align}
\end{lem}
\beg{proof}
By \eqref{ine-b-si0} and It\^o's formula, we have that 
\beg{align}\label{ine1}
\d |X_t^{\mu,x}|^2&=2\<b(X_t^{\mu,x},\mu),X_t\>\d t+2\<X_t^{\mu,x},\si(X_t^{\mu,x},\mu)\d W_t\>+\|\si(X_t^{\mu,x},\mu)\|_{HS}^2\d t\nonumber\\
&\leq \left\{-C_1|X_t^{\mu,x}|^{1+r_1}+C_2+C_3\|\mu\|_{1+r_2}^{r_3}\right\}\d t+2\<X_t^{\mu,x},\si(X_t^{\mu,x},\mu)\d W_t\>.
\end{align}
Then
\beg{align}\label{ine2}
C_1\E\int_0^s |X_t^{\mu,x}|^{1+r_1}\d t+\E |X_s^{\mu,x}|^2\leq |x|^2+\left(C_2+C_3\|\mu\|_{1+r_2}^{r_3}\right)s,~x\geq 0.
\end{align}
By the B-D-G inequality, \eqref{ine1} and \eqref{ine-b-si00}, we have that 
\beg{align*}
\E\sup_{s\in [0,T]}|X_s^{\mu,x}|^2&\leq |x|^2+\left(C_2+C_3\|\mu\|_{1+r_2}^{r_3}\right)T\\
&\qquad+\E\left(\int_0^T|X_t^{\mu,x}|^2\|\si^*(X_t^{\mu,x},\mu)\|^2\d t\right)^{\ff 1 2}\\
&\leq  |x|^2+\left(C_2+C_3\|\mu\|_{1+r_2}^{r_3}\right)T\\
&\qquad +\E\sup_{s\in [0,T]}|X_s^{\mu,x}|\left(\int_0^T \|\si^*(X_t^{\mu,x},\mu)\|^2\d t\right)^{\ff 1 2}\\
&\leq  |x|^2+\left(C_2+C_3\|\mu\|_{1+r_2}^{r_3}\right)T+\ff 1 2 \E\sup_{s\in [0,T]}|X_s^{\mu,x}|^2\\
&\qquad  +\ff 1 2\E\int_0^T\left(C_4|X_t^{\mu,x}|^{1+r_1}+C_5\|\mu\|_{1+r_2}^{r_3}+C_6\right)\d t.
\end{align*}
Combining this with \eqref{ine2}, we obtain \eqref{ine-lem2}.

\end{proof}

Lemma \ref{lem2} together with Lemma \ref{lem3} indicates that $\cT$ is well-defined and  $\sP^{1+r_2}_{M}$ is an invariant subset of $\cT$ for large $M$.
\beg{lem}\label{lem2}
Assume that (H1) and \eqref{ine-C-r}, \eqref{ine-C-r'} and \eqref{inesi} hold. Then for every $\mu\in\sP^{1+r_2}(\R^d)$, \eqref{equ0} has an invariant probability measure. If the invariant probability measure is unique, then there exists $M_0>0$ such that for all $M\geq M_0$, $\cT$ maps $\sP^{1+r_2}_{M}$  into $\sP^{1+r_2}_{M}$. 
\end{lem}

\beg{proof}
Let $q=1+\ff {r_2-r_1} 2$. Then 
\beg{align}\label{eq-qr}
q\geq 1,\qquad q-1=\ff {r_2-r_1} 2,\qquad 2q+r_1-1=1+r_2.
\end{align}
It follows from \eqref{ine1} and \eqref{eq-qr} that
\beg{align*}
\d |X_t^{\mu,x}|^{2q}&\leq q|X_t^{\mu,x}|^{2(q-1)}\left(-C_1|X_t^{\mu,x}|^{1+r_1}+C_2+C_3\|\mu\|_{1+r_2}^{r_3}\right)\d t\\
&\qquad +2q|X_t^{\mu,x}|^{2(q-1)}\<X_t^{\mu,x},\si(X_t^{\mu,x},\mu)\d W_t\>\\
&\leq \left\{-C_1q|X_t^{\mu,x}|^{1+r_2}+q|X_t^{\mu,x}|^{r_2-r_1}\left(C_2+C_3\|\mu\|_{1+r_2}^{r_3}\right)\right\}\d t\\
&\qquad +2q(q-1)|X_t^{\mu,x}|^{r_2-r_1}\left(C_4|X_t^{\mu,x}|^{1+r_2}+C_5\|\mu\|_{1+r_2}^{r_3}+C_6\right)\d t\\
&\qquad +2q|X_t^{\mu,x}|^{2(q-1)}\<X_t^{\mu,x},\si(X_t^{\mu,x},\mu)\d W_t\>\\
&\leq \left\{-C_1q|X_t^{\mu,x}|^{1+r_2}+q|X_t^{\mu,x}|^{r_2-r_1}\left(C_2+C_3\|\mu\|_{1+r_2}^{r_3}\right)\right\}\d t\\
&\qquad +2q(q-1) C_4|X_t^{\mu,x}|^{r_2+1}+2q(q-1) |X_t^{\mu,x}|^{r_2-r_1}\left\{C_5\|\mu\|_{1+r_2}^{r_3}+C_6 \right\}\d t\\
&\qquad +2q|X_t^{\mu,x}|^{r_2-r_1}\<X_t^{\mu,x},\si(X_t^{\mu,x},\mu)\d W_t\>\\
&\leq -q(C_1-2(q-1)C_4)|X_t^{\mu,x}|^{1+r_2}+q(C_2+2C_6(q-1))|X_t^{\mu,x}|^{r_2-r_1} \d t\\
&\qquad + q(C_3+2(q-1)C_5)|X_t^{\mu,x}|^{r_2-r_1} \|\mu\|_{1+r_2}^{r_3} \d t\\
&\qquad +2q|X_t^{\mu,x}|^{r_2-r_1}\<X_t^{\mu,x},\si(X_t^{\mu,x},\mu)\d W_t\>.
\end{align*}
Then
\beg{align}\label{ine-ph0}
&|X_s^{\mu,x}|^{2+r_2-r_1}+q(C_1-(r_2-r_1)C_4)\int_0^s|X_t^{\mu,x}|^{1+r_2}\d t\nonumber\\
&\qquad = |X_s^{\mu,x}|^{2+r_2-r_1}+q(C_1-2(q-1)C_4)\int_0^s|X_t^{\mu,x}|^{1+r_2}\d t\nonumber\\
&\qquad \leq |x|^{2+r_2-r_1}+q(C_2+2C_6(q-1))\int_0^s|X_t^{\mu,x}|^{r_2-r_1}\d t\nonumber\\
&\qquad \qquad +q(C_3+2(q-1)C_5)\|\mu\|_{1+r_2}^{r_3} \int_0^s|X_t^{\mu,x}|^{r_2-r_1}\d t\nonumber\\
&\qquad\qquad +2q\int_0^s|X_t^{\mu,x}|^{2(q-1)}\<X_t^{\mu,x},\si(X_t^{\mu,x},\mu)\d W_t\>,
\end{align}
where we use \eqref{eq-qr} again in the first equality. Since $r_2-r_1<1+r_2$ and \eqref{ine-C-r}, there exist positive constants $\tld C_1,\tld C_2$  such that
\beg{align*}
\E|X_s^{\mu,x}|^{2+r_2-r_1}+\int_0^s|X_t^{\mu,x}|^{r_2+1}\d t\leq \tld C_1 |x|^{2+r_2-r_1}+\tld C_2t,~t\geq 0,
\end{align*}
which yields that
\beg{align}\label{w-con-P}
\ff 1 t\int_0^t\E |X_s^{\mu,x}|^{r_2+1}\d s\leq \tld C_2+\ff {\tld C_1} t  |x|^{2+r_2-r_1},~t>0.
\end{align}
Since $|\cdot|^{1+r_2}$ is a compact function, the existence of the invariant probability measure of \eqref{equ0} follows from the Bogoliubov-Krylov theorem.

If the invariant probability measure is unique, then $\cT_\mu$ is well defined. It follows from \eqref{w-con-P} that   there exists a sequence $ t_n \uparrow +\infty$   such that as $n\ra +\infty$,
\beg{align*}
\ff 1 {t_n}\int_0^{t_n} P_s^{\mu,*}\de_0\d s \xrightarrow{w} \cT_\mu.
\end{align*}
Thus
\beg{align*}
\lim_{n\ra+\infty}\ff 1 {t_n}\int_0^{t_n}  P_s^{\mu,*}\de_0( |\cdot|^{1+r_2}\we N)\d s  = \cT_\mu(|\cdot|^{1+r_2}\we N),~N\geq 1,
\end{align*}
which  together with \eqref{w-con-P} yields that
\beg{align*}
\cT_\mu(|\cdot|^{1+r_2}\we N)\leq \tld C_2,~N\geq 1.
\end{align*}
By Fatou's lemma, $\cT_\mu(|\cdot|^{1+r_2})<\infty$. 
Then, we can integrate both side of  \eqref{ine-ph0} w.r.t. $\cT_\mu$, and get that 
\beg{align*}
&\cT_\mu(|\cdot|^{2+r_2-r_1})+q(C_1- (r_2-r_1)C_4)\int_0^t\cT_\mu(|\cdot|^{1+r_2})\d s\\
&\qquad \leq \cT_\mu(|\cdot|^{2+r_2-r_1})+q(C_3+2(q-1)C_5)\|\mu\|_{1+r_2}^{r_3}\cT_\mu(|\cdot|^{r_2-r_1})t\\
&\qquad \qquad +q(C_2+2C_6(q-1))\cT_\mu(|\cdot|^{r_2-r_1})t.
\end{align*} 
This combining with \eqref{ine-C-r} and \eqref{eq-qr} implies by the Jensen inequality that
\beg{align}\label{ine-F-mu}
\cT_\mu(|\cdot|^{1+r_2})&\leq \ff {C_2+C_6(r_2-r_1)} {C_1-C_4(r_2-r_1)}\|\cT_\mu\|_{1+r_2}^{ {r_2-r_1} }\nonumber\\
&\qquad +\ff {C_3+C_5(r_2-r_1)} {C_1-C_4(r_2-r_1)}\|\cT_\mu\|_{1+r_2}^{ {r_2-r_1} }\|\mu\|_{1+r_2}^{r_3}.
\end{align}
It is clear that $\mu$ with $\cT_\mu(|\cdot|^{1+r_2})=0$ always satisfies \eqref{ine-F-mu}. Thus we can   assume that $\cT_\mu(|\cdot|^{1+r_2})\neq 0$. Then we arrive at
\beg{align}\label{ine-F-mu1}
\|\cT_\mu\|_{1+r_2}^{  {1+r_1} }&\leq \ff {C_2+C_6(r_2-r_1)} {C_1-C_4(r_2-r_1)} +\ff {C_3+C_5(r_2-r_1)} {C_1-C_4(r_2-r_1)}\|\mu\|_{1+r_2}^{r_3}.
\end{align}
If  $r_3= {1+r_1} $, then by  \eqref{ine-C-r'}, for 
$$M\geq \left(\ff {C_2+C_6(r_2-r_1)} {C_1-C_3-(C_4+C_5)(r_2-r_1)}\right)^{\ff {1  } {r_3}}=: M_0,$$
we have that
\beg{align}\label{ine-M}
\ff {C_2+C_6(r_2-r_1)} {C_1-C_4(r_2-r_1)} +\ff {C_3+C_5(r_2-r_1)} {C_1-C_4(r_2-r_1)}M^{r_3}\leq M^{  {1+r_1} }. 
\end{align}
Consequently, $\|\cT_\mu\|_{1+r_2}  \leq M$.\\
If $r_3< 1+r_1 $, then it is easy to see that there exists $M_0$ depending on $C_i$ ($i=1,\cdots, 6$) and $r_1,r_2,r_3$ such that \eqref{ine-M} holds for each $M\geq M_0$. Hence, for  every $\mu\in \sP^{1+r_2}_{M}$,
we have from \eqref{ine-F-mu1} that 
\beg{align*}
&\|\cT_\mu\|_{1+r_2}^{1+r_1}\leq \ff {C_2+C_6(r_2-r_1)} {C_1-C_4(r_2-r_1)} + \ff {C_3+C_5(r_2-r_1)} {C_1-C_4(r_2-r_1)} M^{r_3}\leq M^{ 1+r_1 }.
\end{align*}
Therefore, we prove that  $\cT_\mu\in \sP^{1+r_2}_{M}$ for every $M\geq M_0$.

\end{proof}


We not only prove the uniqueness of the probability measure to $P_t^\mu$, but also establish the so called $V$-uniformly exponential ergodicity.
\beg{defn}
The Markov $P_t^\mu$ is $V$-uniformly exponential ergodic, if there exist $C>0$ and $\ga>0$ such that  
\beg{align*}
\sup_{\|f\|_V\leq 1}\left|P_t^\mu f(x)-\mu(f)\right|\leq CV(x)e^{-\ga t},~x\in\R^d.
\end{align*}
\end{defn}
We have the following lemma. 
\beg{lem}\label{lem3}
Fix $\mu\in\sP^{1+r_2}$. Assume (H1), \eqref{ine-b-si00}, \eqref{ineb-g} and that $\si(\cdot,\mu)$ is continuous and non-degenerate. Then $P_t^\mu$ has at most one  invariant probability measure and  is $V$-uniformly exponential ergodic with 
\beg{equation*}
V(x)=\beg{cases}
e^{\de(1+|X_t|^2)^{\ff {1-r_1} 2}},~&\text{if}~r_1<1;\\
1+|x|,~&\text{if}~r_1=1;\\
1,~&\text{if}~r_1>1.
\end{cases}
\end{equation*}
\end{lem}
\beg{proof}
Under the assumption of this lemma, it follows from \cite[Lemma 7.3]{XZ} that $P_t^\mu$ is the strong Feller property and irreducibility. The $V$-uniformly exponential ergodicity has been proved in \cite[Theorem 2.9]{XZ} for   $r_1\geq 1$. 

For $r_1\in [0,1]$, it follows from (H1) that there exist positive constants $\bar C_1, \bar C_2, \bar C_3$  depending on $\mu(|\cdot|^{r_2+1})$ such that
\beg{align*}
2\<x,b(x,\mu)\>+\|\si(x,\mu)\|_{HS}^2&\leq -\bar C_1(1+|x|^2)^{\ff {r_1+1} 2}+\bar C_2\\
\|\si(x,\mu)\|_{HS}^2&\leq  \bar C_3(1+|x|^2)^{r_1},~x\in\R^d.
\end{align*}
Then  we have by the It\^o formula that 
\beg{align*}
\d \left(1+|X_t^{\mu,x}|^2\right)^{\ff {1-r_1} 2}&\leq -\left(\ff {(1-r_1)\bar C_1} 2-\ff {\bar  C_2(1-r_1)} {2(1+|X_t^{\mu,x}|^2)}+\ff {(1-r_1^2)|\si^*(X_t^{\mu,x})X_t|^2} {2(1+|X_t^{\mu,x}|^2)^{\ff {3+r} 2}}\right)\d t\\
&\qquad +\ff {(1-r_1)\<X_t^{\mu,x},\si(X_t^{\mu,x})\d W_t\>} {(1+|X_t^{\mu,x}|^2)^{\ff {1+r_1} 2}},~x\in\R^d.
\end{align*}
Let  
$$\de=\ff {\bar C_1} {2\bar C_3(1+r_1)},~\hat C_2=\ff {\de \bar C_2(1-r_1)} 2 \sup_{1\leq u\leq (\ff {4\bar C_2} {C_1})^{\ff 1 {1+r_1}}}  u^{-(1+r_1)}e^{\de u^{1-r_1}},$$
and let $V(x)=e^{\de(1+|x|^2)^{\ff {1-r_1} 2}}$. Then
\beg{align*}
&\d V(X_t^{\mu,x})-\ff {\de (1-r_1) V(X_t^{\mu,x})\<X_t^\mu,\si(X_t^{\mu,x})\d W_t\>} {(1+|X_t^{\mu,x}|^2)^{\ff {1+r_1} 2}}\\
&\qquad \leq -\de V(X_t^\mu)\left(\ff {\bar  C_1(1-r_1)} 2+\ff {\bar  C_2(1-r_1)} {2(1+|X_t^{\mu,x}|^2)}-\ff {\de(1-r_1)^2|\si^*(X_t^\mu)X_t^{\mu,x}|^2} {2(1+|X_t^{\mu,x}|^2)^{1+r_1}}\right)\d t\\
&\qquad\leq -\de V(X_t^\mu) \left(\ff {\bar  C_1(1-r_1)} 2-\ff {\de(1-r_1)^2\bar C_3} 2-\ff {\bar C_2(1-r_1)} {2(1+|X_t^{\mu,x}|^2)^{\ff {1+r_1} 2}}\d t\right) \d t\\
&\qquad\leq \de V(X_t^{\mu,x})\left(- \ff {\bar C_1(1-r_1)} 4 +\ff {\bar C_2(1-r_1)} {2(1+|X_t^{\mu,x}|^2)^{\ff {1+r_1} 2}}\1_{[\ff {4\bar C_2} {\bar C_1}\leq (1+|X_t^{\mu,x}|^2)^{\ff {1+r_1} 2} ]} \right)\d t\\
&\qquad \qquad +\sup_{\ff {4\bar C_2} {\bar C_1}\geq (1+|x|^2)^{\ff {1+r_1} 2} }\left(\ff {\bar C_2(1-r_1)V(x)} {2(1+|x|^2)^{\ff {1+r_1} 2}}\right)\d t\\
&\qquad \leq -\ff {\bar C_1(1-r_1)\de} 8 V(X_t^{\mu,x})\d t+\hat C_2\d t.
\end{align*}
Hence, 
\beg{align*}
\E e^{\de(1+|X_t^{\mu,x}|^2)^{\ff {1-r_1} 2}}\leq e^{-\ff {\bar C_1(1-r_1)\de t} 8} e^{\de (1+|x|^2)^{\ff {1-r_1} {2}}}+\ff {8\hat C_2} {\bar  C_1(1-r_1)\de},~t\geq 0,~x\in\R^d. 
\end{align*}
Combining this with the strong Feller property and irreducibility of $P_t^\mu$, we can prove that $P_t^\mu$ is $V$-uniformly exponential ergodicity by following the proof of \cite[Theorem 2.5.]{GoMa} line to line, and we omit the details here. 

\end{proof}

It is clear that $\sP^{1+r_2}_{M}$ is a convex and  weak compact set in $ \sP^1 $. Next, we prove that $\cT$ is weak continuous in $\sP^{1+r_2}_{M}$.

\beg{lem}\label{lem4}
Under the assumption of Theorem \ref{thm1},  $\cT$ is weak continuous in $\sP^{1+r_2}_{M}$.
\end{lem}
\beg{proof}
For all $f\in \sB_b(\R^d)$ and $t>0$, we have that  
\beg{align}\label{TT}
\left|\cT_\mu(f)-\cT_\nu(f)\right|&=\left|\cT_\mu(f)-\cT_\nu(P_t^\nu f)\right|\nonumber\\
&\leq \left|\cT_\nu(\cT_\mu(f))-\cT_\nu(P_t^\mu f)\right|+\left|\cT_\nu(P_t^\mu f)- \cT_\nu(P_t^\nu f)\right|\nonumber\\
&\leq \cT_\mu\left(\left|f- P_t^\mu f \right|\right)+\cT_\nu\left(\left| P_t^\mu f- P_t^\nu f \right|\right).
\end{align}

By Lemma \ref{lem3}, there exist positive constants $\ga$ and $C$ depending on $\mu$ such that
\beg{align*}
\left|\cT_\mu(f)-P_t^\mu f(x)\right|\leq CV(x) e^{-\ga t},~x\in\R^d.
\end{align*}
Then for each $m>0$, we have by Lemma \ref{lem2} and Lemma \ref{lem3} that 
\beg{align*}
\cT_\nu\left(\left|f- P_t^\mu f \right|\right)&\leq \cT_\nu\left(CV(\cdot)e^{-\ga t}\1_{[|\cdot|\leq m]}+\|f\|_\infty\1_{[|\cdot|> m]}\right)\\
&\leq C\left(\sup_{|x|\leq m}V(x)\right)e^{-\ga t}+\|f\|_\infty\cT_\nu(|\cdot|>m)\\
&\leq C\left(\sup_{|x|\leq m}V(x)\right)e^{-\ga t}+\ff {\|f\|_\infty} {m }\|\cT_\nu\|_{1+r_2}\\
&\leq C\left(\sup_{|x|\leq m}V(x)\right)e^{-\ga t}+\ff {\|f\|_\infty M} {m }. 
\end{align*}
Hence, 
\beg{align*}
\varlimsup_{t\ra+\infty}\sup_{\nu\in \sP^{1+r_2}_{M}}\cT_\nu\left(\left|f- P_t^\mu f \right|\right)=0.
\end{align*}
Then, for each $\ep>0$, we choose  $t_\ep>0$ such that  
\beg{align}\label{Tf-Pf}
\sup_{\nu\in \sP^{1+r_2}_{M}}\cT_\nu\left(\left|f- P_{t_\ep}^\mu f \right|\right)<\ep.
\end{align}

Let 
$$\ta_n =\inf\{~t>0~|~|X_t^{\mu,x}|\vee |X_t^{\nu,x}|\geq n\}.$$
For every $f\in C_b(\R^d)\cap Lip(\R^d)$,  we have that
\beg{align*}
\cT_\nu\left(\left| P_{t_\ep}^\mu f- P_{t_\ep}^\nu f \right|\right)&\leq \int_{\R^d}\E\left|f(X_{t_\ep}^{\mu,x})-f(X_{t_\ep}^{\nu,x})\right|\cT_\nu(\d x)\\
&\leq \int_{\R^d}\E\left(\left|f(X_{t_\ep}^{\mu,x})-f(X_{t_\ep}^{\nu,x})\right|\1_{[  t_\ep<\ta_n ] }\right)\cT_\nu(\d x)\\
&\qquad +\int_{\R^d}\E\left(\left|f(X_{t_\ep}^{\mu,x})-f(X_{t_\ep}^{\nu,x})\right|\1_{  [t_\ep\geq \ta_n ]  }\right)\cT_\nu(\d x)\\
&\leq \|\nn f\|_\infty\int_{|x|\leq \ff n 2}\E\left(\left|X_{t_\ep}^{\mu,x}-X_{t_\ep}^{\nu,x}\right|\1_{[t_\ep<\ta_n]}\right)\cT_\nu(\d x)\\
&\qquad +2\|f\|_\infty\cT_\nu(|\cdot|\geq \ff  n 2)+2\|f\|_\infty\cT_\nu\left(\P(t_\ep\geq \ta_n )\right)\\
&\leq \|\nn f\|_\infty\int_{|x|\leq \ff n 2}\E\left(\left|X_{t_\ep\we\ta_n}^{\mu,x}-X_{t_\ep\we\ta_n}^{\nu,x}\right|\right)\cT_\nu(\d x)\\
&\qquad +\ff 4 {n}\|f\|_\infty \|\cT_\nu\|_{1+r_2}\\
&\qquad +2\|f\|_\infty\int_{\R^d} \P\left(\sup_{s\in[0,t_\ep]}|X_s^{\mu,x}|\vee|X_s^{\nu,x}| \geq n\right)  \cT_\nu(\d x)\\
&=: I_{1,n}+I_{2,n}+I_{3,n}.
\end{align*}

For $I_{2,n}$. Due to Lemma \ref{lem2} and $\nu\in \sP^{1+r_2}_{M}$, we have that $\|\cT_\nu\|_{1+r_2}\leq M$. Then 
\beg{align*}
\lim_{n\ra+\infty}\sup_{\nu\in\sP^{1+r_2}_{M}}I_{2,n}=0.
\end{align*}

For $I_{3,n}$. By Lemma \ref{lem1}, the Chebyshev inequality and that $\nu,\mu\in\sP^{1+r_2}_{M}$, there exist $C>0$ depending on $C_i$, $i=1,\cdots, 6$, $r_1,r_2,r_3$, $M$ and $t_\ep$ such that
\beg{align*}
&\P\left(\sup_{s\in[0,t_\ep]}|X_s^{\mu,x}|\vee|X_s^{\nu,x}| \geq n\right)\\
&\qquad \leq \left(\P\left(\sup_{s\in[0,t_\ep]}|X_s^{\mu,x}|\geq n\right)+\P\left(\sup_{s\in[0,t_\ep]}|X_s^{\nu,x}|\geq n\right)\right)\\
&\qquad \leq n^{-(2\we (1+r_2))}\E\left(\sup_{s\in[0,t_\ep]}|X_s^{\mu,x}|^{2\we(1+r_2)}+ \sup_{s\in[0,t_\ep]}|X_s^{\nu,x}|^{2\we(1+r_2)}\right)\\
&\qquad \leq C(1+|x|^{2\we(1+r_2)})n^{-(2\we (1+r_2))}.
\end{align*}
Since $\cT_\nu\in \sP^{1+r_2}_{M}$, it is clear that
\beg{align*}
\lim_{n\ra+\infty}\sup_{\nu\in\sP^{1+r_2}_{M}}\cT_\nu(I_{3,n})=0.
\end{align*}

For $I_{1,n}$. It is easy to see that
\beg{align}\label{I2n}
\left|X_s^{\mu,x}-X_s^{\nu,x}\right|&\leq \int_0^s\left|b(X_t^{\mu,x},\mu)-b(X_t^{\nu,x},\nu)\right|\d t+\left|\int_0^s\left(\si(X_t^{\mu,x},\mu)-\si(X_t^{\nu,x},\nu)\right)\d W_t\right|\nonumber\\
&\leq \int_0^s\left|b(X_t^{\mu,x},\mu)-b(X_t^{\nu,x},\mu)\right|\d t+\int_0^s\left|b(X_t^{\nu,x},\mu)-b(X_t^{\nu,x},\nu)\right|\d t\nonumber\\
&\qquad +\left|\int_0^s\left(\si(X_t^{\mu,x},\mu)-\si(X_t^{\nu,x},\mu)\right)\d W_t\right|\nonumber\\
&\qquad +\left|\int_0^s\left(\si(X_t^{\nu,x},\mu)-\si(X_t^{\nu,x},\nu)\right)\d W_t\right|,~s<\ta_n.
\end{align}
It follows from (H2)  that
\beg{align*}
\int_0^{s\we\ta_n}\left|b(X_t^{\mu,x},\mu)-b(X_t^{\nu,x},\mu)\right|\d t\leq K_n \int_0^{s\we\ta_n}\left|X_t^{\mu,x} - X_t^{\nu,x}\right|\d t
\end{align*}
and 
\beg{align*}
&\E\sup_{r\in[0,s]}\left|\int_0^{r\we\ta_n}\left(\si(X_t^{\mu,x},\mu)-\si(X_t^{\nu,x},\mu)\right)\d W_t\right|\\
&\qquad \leq \E \left(\int_0^{s\we\ta_n}\left\|\si(X_t^{\mu,x},\mu)-\si(X_t^{\nu,x},\mu) \right\|_{HS}^2\d t\right)^{\ff 1 2}\\
&\qquad \leq K_n \E \left(\int_0^{s\we\ta_n}\left| X_t^{\mu,x} -X_t^{\nu,x}  \right|^2\d t\right)^{\ff 1 2}\\
&\qquad \leq K_n \E\left(\sup_{t\in [0,s\we\ta_n]}\left| X_{t }^{\mu,x} -X_{t }^{\nu,x}  \right|^{\ff 1 2} \right)\left(\int_0^{s\we\ta_n}\left| X_t^{\mu,x} -X_t^{\nu,x}  \right|\d t\right)^{\ff 1 2}\\
&\qquad \leq \ff 1 2\E \sup_{t\in [0,s ]}\left| X_{t\we\ta_n }^{\mu,x} -X_{t\we\ta_n}^{\nu,x}  \right|  +\ff 1 2K_n^2 \E \int_0^{s\we\ta_n}\left| X_t^{\mu,x} -X_t^{\nu,x}  \right|\d t. 
\end{align*}
Putting these into \eqref{I2n}, we have by using the Gronwall inequality that
\beg{align}\label{I1n1}
&\E \sup_{t\in [0,s]}\left| X_{t\we\ta_n }^{\mu,x} -X_{t\we\ta_n}^{\nu,x}  \right|\nonumber\\
&\qquad  \leq 2e^{(2K_n +K_n^2)s}\Big(\E\sup_{r\in [0,s]} \int_0^{r\we\ta_n}\left|b(X_t^{\nu,x},\mu)-b(X_t^{\nu,x},\nu)\right|\d t\nonumber\\
&\qquad\qquad  +\E\sup_{r\in [0,s]} \left|\int_0^{r\we\ta_n}\left(\si(X_t^{\nu,x},\mu)-\si(X_t^{\nu,x},\nu)\right)\d W_t\right|\Big),~s>0. 
\end{align}
Because 
\beg{align*}
\E\sup_{r\in [0,s]} \int_0^{r\we\ta_n}\left|b(X_t^{\nu,x},\mu)-b(X_t^{\nu,x},\nu)\right|\d t&\leq \E\int_0^{s}\left|b(X_{t\we\ta_n}^{\nu,x},\mu)-b(X_{t\we\ta_n}^{\nu,x},\nu)\right|\d t\\
&\leq s\sup_{|x|\leq n}|b(x,\mu)-b(x,\nu)| 
\end{align*}
and 
\beg{align*}
&\E\sup_{r\in [0,s]} \left|\int_0^{r\we\ta_n}\left(\si(X_t^{\nu,x},\mu)-\si(X_t^{\nu,x},\nu)\right)\d W_t\right|\\
&\qquad \leq \E\left(\int_0^{s\we\ta_n}\|\si(X_t^{\nu,x},\mu)-\si(X_t^{\nu,x},\nu)\|_{HS}^2\d t\right)^{\ff 1 2}\\
&\qquad \leq \E\left(\int_0^{s }\|\si(X_{t\we\ta_n}^{\nu,x},\mu)-\si(X_{t\we\ta_n}^{\nu,x},\nu)\|_{HS}^2\d t\right)^{\ff 1 2} \\
&\qquad \leq \sq s\sup_{|x|\leq n}\|\si(x,\mu)-\si(x,\nu)\|_{HS},
\end{align*}
we have that
\beg{align*}
I_{1,n}&\leq \|\nn f\|_\infty\int_{|x|\leq \ff n 2}\E \sup_{t\in [0,t_\ep]}\left| X_{t\we\ta_n }^{\mu,x} -X_{t\we\ta_n}^{\nu,x}  \right|\cT_\nu(\d x)\\
&\leq 2e^{(2K_n(M)+K_n(M)^2)t_\ep} \|\nn f\|_\infty\\
&\qquad \times\left(t_\ep\sup_{|x|\leq n}|b(x,\mu)-b(x,\nu)| +\sq{t_\ep}\sup_{|x|\leq n}\|\si(x,\mu)-\si(x,\nu)\|_{HS}\right).
\end{align*}
This together with \eqref{ineb} implies that for $\nu\xrightarrow{w}\mu$ in $\sP^{1+r_2}_{M}$
\beg{align*}
\lim_{\nu\xrightarrow{w}\mu}I_{1,n}=0.
\end{align*}

Letting $\nu\xrightarrow{w} \mu$ in $\sP^{1+r_2}_{M}$ and $n\ra+\infty$, we get that 
\beg{align*}
\lim_{\nu\xrightarrow{w}\mu}\cT_\nu\left(\left| P_{t_\ep}^\mu f- P_{t_\ep}^\nu f \right|\right)=0.
\end{align*}
This together with \eqref{Tf-Pf} and \eqref{TT} yields that for all $\ep>0$ 
\beg{align*}
\lim_{\nu\xrightarrow{w}\mu}\left|\cT_\mu(f)-\cT_\nu(f)\right|\leq \ep.
\end{align*}
Therefore, the proof is completed. 

\end{proof}

\noindent {\bf{Proof of Theorem \ref{thm1}.~~}} By Lemma \ref{lem2}, there is some $M_0>0$ such that for each $M\geq M_0$, $\sP^{1+r_2}_{M}$ is   invariant under $\cT$.   By \cite[Theorem 5.5]{Chen} or \cite[Theorem 5.5]{CarDe}, $\sP^{1+r_2}_{M}$ is a convex and compact subset of $(\sP^1,W)$. By Lemma \ref{lem4}, $\cT$ is weak continuous in $\sP^{1+r_2}_{M}$. Since $\sP^{1+r_2}_{M}$ is  compact in $\sP_1$, $\cT$ is continuous in $\sP^{1+r_2}_{M}$ w.r.t. to the Kantorovich-Rubinstein norm.  Hence,  the Schauder fixed point theorem yields that $\cT$ has a fixed in $\sP^{1+r_2}_{M}$.

\qed




\bigskip

\noindent{\bf Proof of Corollary \ref{cor0}}

Because $\si$ is bounded and satisfies (H2) and (H3), we  focus on the drift term below when verifying (H1)-(H3).  It follows from \eqref{VF0} that
\beg{align}\label{nn-VF0}
&\<\nn V(x)-\nn V(0),x\>+\mu(\<\nn F(x-\cdot)-\nn F(0-\cdot),x\>) \\
&\quad = \int_0^1\int_{\R^d}\left(\<\nn^2 V(\th x)x,x\>+ \<\nn^2 F(\th x-y) x,x\>\right)\mu(\d y)\d \th\\
&\quad \geq  \left(\be_2 |x|^2-\be_1|x|+\be_0-\be_3\mu(|\cdot|^2) \right) |x|^2. 
\end{align}
Let $\de_1=\be_3/\be_2$, it follows by the H\"older inequality that 
\beg{align}\label{mu-44}
 \be_3 \mu(|\cdot|^2)  |x|^2\leq \de_1 \be_2|x|^4+\ff {\be_3^2} {4\de_1\be_2}\|\mu\|_4^4.
\end{align}
By \eqref{nnF1} and the H\"older inequality, for  $\de_2=\al_2/(4\be_2)$, we have that    
\beg{align}\label{nnF0x}
|\mu(\<\nn F(0-\cdot),x\>)|&\leq \al_1|x|+\al_2\|\mu\|_3^3|x|\leq \al_1|x|+\al_2\|\mu\|_4^3|x| \nonumber\\
&\leq \al_1|x|+\ff {3\al_2^{\ff 4 3}} {4(4\de_2\be_2)^{\ff 1 3}}\|\mu\|_4^4+\de_2\be_2|x|^4
\end{align}
Putting \eqref{mu-44} and \eqref{mu-44} into \eqref{nn-VF0}, we have that
\beg{align}\label{nn-VF1}
\<\nn V(x)+\nn F*\mu(x),x\>&\geq (1-\de_1-\de_2)\be_2|x|^4 -\be_1|x|^3+\be_0|x|^2-(|\nn V(0)|+\al_1)|x|\nonumber\\
&\qquad -\left(\ff {3\al_2^{\ff 4 3}} {4(4\de_2\be_2)^{\ff 1 3}}+\ff {\be_3^2} {4\de_1\be_2}\right)\|\mu\|_4^4\nonumber\\
&= \left( \be_2-\ff {\al_2}4-\ff {\be_3} 2\right)|x|^4  -\be_1|x|^3+\be_0|x|^2\nonumber\\
&\qquad -(|\nn V(0)|+\al_1)|x|-\left(\ff {3\al_2} {4 }+\ff {\be_3} {2 }\right)\|\mu\|_4^4.
\end{align}
Since $\be_2>\al_2+\be_3$, it holds that 
\beg{align*}
\be_2-\ff {\al_2}4-\ff {\be_3} 2>\ff {3\al_2} {4 }+\ff {\be_3} {2 }.
\end{align*}
This together with \eqref{nn-VF1} and  the H\"older inequality yields that there are positive constants $  C_1, C_2, C_3$ with $\tld C_1>\tld C_3$ such that
\beg{align*}
\<b(x,\mu),x\>& \leq - C_1|x|^4+ C_2+  C_3\|\mu\|_4^4.
\end{align*} 
Combining this with that $\si$ is bounded on $\R^d\times\sP(\R^d)$, we have that  (H1) holds with $r_1=r_2=3,r_3=4$. Moreover, we also have that $C_4=C_5=0$ in \eqref{ine-b-si00} due to that  $\si$ is bounded.  Furthermore, \eqref{ine-C-r} and \eqref{ine-C-r'} hold since $r_1=r_2$,$r_3=1+r_1$, $C_4=C_5=0$ and $ C_1>C_3$.

Because that $V$  is  twice continuous differentiable, $\nn V$ is local Lipschitz. By \eqref{n2F4} and the mean value theorem, we have that 
\beg{align*}
&\mu\left(|\nn F(x-\cdot)-\nn F(y-\cdot)|\right)\\
&\qquad\leq |x-y|\int_0^1\mu\left(\|\nn^2 F(x+\th (y-x)-\cdot)\|\right)\d \th\\
&\qquad\leq \al_3|x-y|\left(1+ 2^{(\ga_0-1)^+}\left(\int_0^1((1-\th)  |x|+\th |y|)^{\ga_0}\d \th +\|\mu\|_{\ga_0}^{\ga_0}\right)\right)\\
&\qquad\leq C_{\al_3,\ga_0}|x-y|\left(1+  |x|^{\ga_0}+ |y|^{\ga_0}  +\|\mu\|_{4}^{\ga_0} \right).
\end{align*}
Thus \eqref{b0-K} holds. The inequality \eqref{ineb-g} follows from \eqref{nnV1} and  \eqref{nnF1} directly.  Hence, (H2) holds.

For any $\pi\in\sC(\mu,\nu)$, we have by \eqref{n2F4} that 
\beg{align*}
|b(x,\mu)-b(x,\nu)| & =\left| \mu(\nn F(x-\cdot))-\nu(\nn F(x-\cdot))\right|\\
&= \left|\int_{\R^d\times\R^d}\left(\nn F(x-y_1)-\nn F(x-y_2)\right)\pi(\d y_1,\d y_2)\right| \\
&\leq  C_{\al_3,\ga_0} \int_{\R^d\times\R^d} \left(1+ |x|^{\ga_0}+  |y_1|^{\ga_0}+ |y_2|^{\ga_0}\right)  |y_1-y_2|\pi(\d y_1,\d y_2)\\
&\leq C_{\al_3,\ga_0}(1+|x|^{\ga_0}) \int_{\R^d\times\R^d} |y_1-y_2| \pi(\d y_1,\d y_2)\\
&\quad +C_{\al_3,\ga_0}\int_{\R^d\times\R^d}(|y_1|^{\ga_0}+|y_2|^{\ga_0})|y_1-y_2|\pi(\d y_1,\d y_2)\\
&\leq   C_{\al_3,\ga_0}(1+|x|^{\ga_0})\int_{\R^d\times\R^d} |y_1-y_2| \pi(\d y_1,\d y_2)\\
&\quad +\tld C_{\al_3,\ga_0} (\|\mu\|_4^{\ga_0}+\|\nu\|_4^{\ga_0})\left(\int_{\R^d\times\R^d} |y_1-y_2|^{\ff 4 {4-\ga_0}}\pi(\d y_1,\d y_2)\right)^{\ff {4-\ga_0} 4}.
\end{align*}
Thus 
\beg{align*}
|b(x,\mu)-b(x,\nu)| \leq C_{\al_3,\ga_0}(1+|x|^{\ga_0})W(\mu,\nu)+\tld C_{\al_3,\ga_0} (\|\mu\|_4^{\ga_0}+\|\nu\|_4^{\ga_0})W_{\ff 4 {4-\ga_0}}(\mu,\nu),
\end{align*}
where $W_{\ff 4 {4-\ga_0} }(\mu,\nu)$ is the $\ff 4 {4-\ga_0} $-Wasserstein distance.  Note that $\ff 4 {4-\ga_0}<1+r_2$ since $\ga_0<3$ and $r_2=3$. For $\mu_m,\mu\in\sP^{1+r_2}_M$ with $\mu_m\xrightarrow{w}\mu$, we have by \cite[Theorem 5.6]{Chen} or \cite[Theorem 5.5]{CarDe} that 
\beg{align*}
\lim_{m\ra+\infty}\left(W(\mu_m,\mu)+W_{\ff 4 {4-\ga_0}}(\mu,\nu)\right)=0.
\end{align*}
Hence, (H3) holds. 

\qed

\subsection{Proof of Theorem \ref{thmN} and Corollary \ref{cor2}}

\noindent{\bf Proof of Theorem \ref{thmN}}

Let 
$$M_0=\sup_{M\geq 0}\left\{~M>0~\Big|~\ff {C_2+C_6(r_2-r_1)} {C_1-C_4(r_2-r_1)} +\ff {C_3+C_5(r_2-r_1)} {C_1-C_4(r_2-r_1)}M^{r_3}\leq M^{ 1+r_1 }\right\}.$$
According to Theorem \ref{thm1}, we find a fixed point of $\cT$ in $\sP^{1+r_2}_{M}$ for every $M\geq M_0$. Let $\sP_{a,\ka }^{M,\ga }=\sP^{1+r_2}_{M}\cap\sP^{\ga_1,\ga_2}_{a,\ka_1,\ka_2}$. We  find the stationary probability in  $\sP_{a,\ka }^{M,\ga }$.

By It\^o's formula, we have that
\beg{align*}
&|X_t^{\mu,a}-a|^2-2\int_0^t \<X_s^{\mu,a}-a,\si(X_s^{\mu,a},\mu)\d W_s\>\\
&\qquad =\int_0^t\left(2\<b_0(X_s^{\mu,a}),X_s^{\mu,a}-a\>+2\<X_s^{\mu,a}-a,b_1(X_s^{\mu,a},\mu)\>\right)\d s\\
&\qquad \leq -\int_0^t g(|X_s^{\mu,a}-a|^\ga_1,\|\mu_a\|_{\ga_1},\|\mu_a\|_{\ga_2})\d s\\
&\qquad = -\int_0^t g_1(|X_s^{\mu,a}-a|^{\ga_1},\|\mu_a\|_{\ga_1},\|\mu_a\|_{\ga_2})\d s\\
&\qquad\qquad -\int_0^t g_2(|X_s^{\mu,a}-a|^{\ga_2},\|\mu_a\|_{\ga_1},\|\mu_a\|_{\ga_2}) \d s.
\end{align*}
Then by Jensen's inequality, we have  that
\beg{align}\label{ineXg0}
\ff {\E|X_t^{\mu,a}-a|^2} t \leq -g\left(\ff 1 t\int_0^t \E|X_s^{\mu,a}-a|^{\ga_1}\d s,\|\mu_a\|_{\ga_1},\|\mu_a\|_{\ga_2}\right)
\end{align}
and
\beg{align}\label{ineXg}
\ff {\E|X_t^{\mu,a}-a|^2} t&\leq -\ff 1 t\int_0^t \E \left(g_1(|X_s^{\mu,a}-a|^{\ga_1},\mu)+g_2(|X_s^{\mu,a}-a|^{\ga_2},\mu)\right)\d s \nonumber\\
&\leq -g_1\left(\ff 1 t\int_0^t \E|X_s^{\mu,a}-a|^{\ga_1}\d s,\|\mu_a\|_{\ga_1},\|\mu_a\|_{\ga_2}\right)\nonumber\\
&\qquad -g_2\left(\ff 1 t\int_0^t \E|X_s^{\mu,a}-a|^{\ga_2}\d s,\|\mu_a\|_{\ga_1},\|\mu_a\|_{\ga_2}\right).
\end{align}
By \eqref{w-con-P} and $\ga_1\leq \ga_2<1+r_2$,  there is a sequence $t_n\uparrow +\infty$ such that 
\beg{align*}
\lim_{n\ra+\infty} \ff 1 {t_n}\int_0^{t_n}\E|X_s^{\mu,a}-a|^{\ga_i}\d s=\cT_\mu(|\cdot-a|^{\ga_i}),~i=1,2.
\end{align*}
By the continuity of $g_i(\cdot,\mu)$, $i=1,2$, we take $t=t_n$  in \eqref{ineXg}, and let $n\ra+\infty$. Then 
\beg{align}\label{g00}
g((\|(\cT_\mu)_a\|_{\ga_1}^{\ga_1},&\|\mu_a\|_{\ga_1},\|\mu_a\|_{\ga_2})\leq 0,\\
g_1(\|(\cT_\mu)_a\|_{\ga_1}^{\ga_1},\|\mu_a\|_{\ga_1},\|\mu_a\|_{\ga_2}) & +g_2(\|(\cT_\mu)_a\|_{\ga_2}^{\ga_2},\|\mu_a\|_{\ga_1},\|\mu_a\|_{\ga_2})\leq 0.\label{g12}
\end{align}
Combining \eqref{g00} with \eqref{sset1}, for every  $\mu\in \sP_{a,\ka }^{M,\ga }$, we have  that $\|(\cT_\mu)_a\|_{\ga_1} \leq\ka_1$. This together with \eqref{g12} and \eqref{sset2}, we have that $\|(\cT_\mu)_a\|_{\ga_2} \leq\ka_2$. Hence, $\sP_{a,\ka }^{M,\ga }$ is an invariant set of the mapping $\cT$.  It is clear that $\sP_{a,\ka }^{M,\ga }$ is also a compact, convex  subset of $(\sP^1,W)$.  By the Schauder fixed point theorem, there exists $\mu\in\sP_{a,\ka }^{M,\ga }$ such that $\mu=\cT_\mu$. 

Let $\mu_1\in\sP_{a_1,\ka }^{M,\ga },\mu_2\in \sP_{a_2,\ka }^{M,\ga }$ be two fixed point of $\cT$. Since $ \ka_1 <  \ff  {|a_1-a_2|}   4 $,   the Chebyshev inequality yields that 
\beg{align*}
\mu_1\left(|\cdot-a_1|\geq \ff {|a_1-a_2|} 2\right)\leq \ff {2 (\mu_1(|\cdot-a_1|^{\ga_1}))^{\ff 1 {\ga_1}}} {|a_1-a_2| }\leq  \ff {2 \ka_1} {|a_1-a_2| }<\ff 1 2,\\
\mu_2\left(|\cdot-a_2|\geq \ff {|a_1-a_2|} 2\right)\leq \ff {2 (\mu_2(|\cdot-a_2|^{\ga_1}))^{\ff 1 {\ga_1}}} {|a_1-a_2| }\leq \ff {2 \ka_1} {|a_1-a_2| }<\ff 1 2.
\end{align*}
Therefore $\mu_1\neq \mu_2$.

\qed

\bigskip

\noindent{\bf Proof of Corollary \ref{cor2}}

It follows from \eqref{VF} that
\beg{align}\label{nn2VF0}
\nn^2 V(x)+\nn^2 F(x-y)&\geq \be_0-2\be_1|x-a|+3\be_2|x-a|^2 \nonumber\\
&\geq \ff 3 2\be_2|x|^2-2\be_1|x|-\left(3\be_2|a|^2+ 2\be_1 |a|  -\be_0 \right).
\end{align} 
By \eqref{LipF}, for any $\ep>0$,
\beg{align}\label{nnF3}
|\nn F(x)|\leq \left(\al_1+\ff {2\al_2^{\ff 3 2}} {3(3\ep)^\ff 1 2}\right)+\ep|x|^3.
\end{align}
By \eqref{nnV1}, \eqref{nn2VF0}, \eqref{nnF3} with $\ep<\ff {3\be_2} 2$ and \eqref{n2F4} with $\ga_0=1$, it follows from  Corollary \ref{cor0} that \eqref{equ00} has stationary distributions  in $\sP^{4}$.

We verify \eqref{bgg1} and \eqref{sset3} for $\ka\in (\ka_0,\ka_1)\cup (\ka_2,+\infty)$ and
\beg{align*}
g(r,r_1)=2\left(\be_2 r^4-\be_1 r^3+ \be_0 r^2-(\al_1+\al_2 r_1)r \right)-\|\si\|_\infty,~r,r_1\geq 0.
\end{align*}
Since $a$ is  the  critical point of $V$, $\nn V(a)=0$. Then  
\beg{align*}
\<\nn V(x+a),x\>=\<\nn V(x+a)-\nn V(a),x\>=\int_0^1 \<\nn^2 V(a+\th x) x,x\>\d \th. 
\end{align*}
Thus 
\beg{align*}
&\<\nn V(x+a),x\>+\mu(\<\nn F(x+a-\cdot),x\>)\\
&\qquad =\int_0^1 \<\nn^2 V(a+\th x) x,x\>\d \th+\mu(\<\nn F(a-\cdot),x\>)\\
&\qquad \qquad +\mu(\<\nn F(x+a-\cdot)-\nn F(a-\cdot),x\>)\\
&\qquad =  \int_0^1\int_{\R^d} \<\left(\nn^2 V(a+\th x)+\nn^2 F(a-y+\th x) \right)x,x\>\mu(\d y)\d \th\\
&\qquad\qquad +\mu(\<\nn F(a-\cdot),x\>)\\
&\qquad \geq \left(\be_0-\be_1|x|+\be_2|x|^2 \right)|x|^2-(\al_1+\al_2\mu(|a-\cdot|))|x|.
\end{align*}
Then
\beg{align*}
2\<b(x+a,\mu),x\>+\|\si(x+a,\mu)\|_{HS}\leq -g(|x|,\|\mu_a\|_1),~x\in\R^d,\mu\in\sP^1.
\end{align*}
It is easy to see that 
\beg{align*}
(\pp_r^2 g)(r,r_1)=2\left(12\be_2 r^2-6\be_1 r+2 \be_0  \right).
\end{align*}
Since $\be_1^2<\ff {8} 3 \be_0\be_2$, we have that $(\pp_r^2 g)(r,r_1)>0$. Thus $g(\cdot,r_1)$ is convex on $[0,+\infty)$. Consider the set 
$$A_{r_1}=\{r\geq 0~|~g(r,r_1)\leq 0\}.$$
For any $\ka\geq 0$, we have that 
\beg{align*}
A_{r_1}& \subset A_{\ka} \\
&=\{r\geq 0~| ~2\left(\be_2 r^4-\be_1 r^3+ \be_0 r^2 \right)\leq  2(\al_1+\al_2 \ka)r+\|\si\|_\infty^2\},~0\leq r_1\leq \ka.
\end{align*}
Since the polynomial $\be_2 r^4-\be_1 r^3+ \be_0 r^2$ is convex, there is $r_0\geq 0$ such that
$$\{r\geq 0~| ~2\left(\be_2 r^4-\be_1 r^3+ \be_0 r^2 \right)\leq  2(\al_1+\al_2 \ka)r+\|\si\|_\infty^2\}=[0,r_0].$$
Since $\be_0>\al_2$, for  $\ka\in (\ka_0,\ka_1)\cup (\ka_2,+\infty)$, we have that 
$$\be_2 \ka^3-\be_1 \ka^2+ (\be_0-\al_2) \ka - \al_1>0.$$
By \eqref{si-be-al}, we have that
\beg{align*}
g(\ka,\ka)&=2\left(\be_2 \ka^4-\be_1 \ka^3+ \be_0 \ka^2-(\al_1+\al_2 \ka)\ka \right)-\|\si\|_\infty^2\\
&=2\left(\be_2 \ka^3-\be_1 \ka^2+ (\be_0-\al_2) \ka - \al_1  \right)\ka-\|\si\|_\infty^2\\
& > 0.
\end{align*}
Thus $r_0<\ka$. This yields that $\cup_{0\leq r_1\leq \ka }A_{r_1}\subset [0,\ka)$. Consequently, \eqref{sset3} holds.

\qed

\subsection{Proofs of examples}

\noindent{\bf Proof of Example \ref{ex-0}~}

Setting  
$$b(x,\mu)=-\be(x-a_1)x(x-a_2)-\al \int_{\R}(x-y)\mu(\d y),$$
it follows from Corollary \ref{cor2} that   \eqref{eq-ex1} has  stationary distributions. 

Next, we use Corollary \ref{cor1} to prove the existence of two distinct stationary distributions. Let 
$$g (r,r_1)=2\be(r^4-|2a_1-a_2|r^3+(a_1(a_1-a_2)+\ff {\al} {\be})r^2-\ff {\al r} {\be}r_1-\|\si\|_\infty^2,~r\geq 0.$$
Then for given $r_1\geq 0$
$$(\pp_r^2 g)(r,r_1)=4\be\left(6r^2-3|2a_1-a_2|r+(a_1(a_1-a_2)+ \ff {\al} {\be})\right).$$
By \eqref{aaba}, we have that  
\beg{align*}
3(2a_1-a_2)^2-8 a_1(a_1-a_2)\leq a_1^2+a_2^2+2(a_1-a_2)^2+a_1^2\vee a_2^2<\ff {8\al} {\be}.
\end{align*}
This implies that  $(\pp_r^2 g)(r,r_1)>0$. Thus $g (\cdot,r_1)$ is convex. Moreover,  
\beg{align}\label{bs-g}
2\<b&(x+a_1,\mu),x\>+\|\si(x,\mu)\|_{HS}^2\nonumber\\
&= -2\be x^2(x+a_1)(x+a_1-a_2)-2\al x\int_{\R}\left( x+a_1 -y\right)\mu(\d y) +\|\si\|_\infty^2\nonumber\\
&= -2\be\left(x^4+(2a_1-a_2)x^3+(a_1(a_1-a_2)+\ff {\al} {\be})x^2-\ff {\al x} {\be}\int_\R(a_1-y)\mu(\d y)\right)+\|\si\|_\infty^2\nonumber\\
&\leq -g (|x|,\|\mu_{a_1}\|_1).
\end{align}  
Let  
\beg{align}\label{pp}
p(r) =2\be(r^4- |2a_1-a_2|r^3+(a_1(a_1-a_2)+\ff {\al} {\be})r^2).
\end{align}
Then $p$ is is convex on $[0,+\infty)$. Thus for $\ka>0$, there is $r_0>0$ such that 
$$\left\{r\geq 0~|~p(r)\leq  2\al\ka r   +\|\si\|_\infty^2 \right\}=[0,r_0].$$
Hence, for any $M>0$, $\mu\in\sP^{ M,1}_{a_1,\ka}$ and $0\leq r_1\leq\ka$, we have that
\beg{align}\label{Ar0}
A_{r_1} &:=\left\{r\geq 0~|~g(r,r_1)\leq 0 \right\}\nonumber\\
&\subset \left\{r\geq 0~|~2\be (r^4-|2a_1-a_2|r^3+(a_1(a_1-a_2)+\ff {\al} {\be})r^2-\ff {\al\ka r} {\be})-\|\si\|_\infty^2\leq 0\right\}\nonumber\\
&=\left\{r\geq 0~|~p(r)\leq  2\al\ka r   +\|\si\|_\infty^2 \right\}\nonumber\\
&=[0,r_0].
\end{align}
Notice that
\beg{align*}
p(\ka)-  2\al  \ka^2-\|\si\|_\infty^2&=2\be(\ka^2- |2a_1-a_2|\ka +a_1(a_1-a_2) )\ka^2-\|\si\|_\infty^2.
\end{align*}
Since $a_1a_2<0$, we have that $a_1(a_1-a_2)>0$ and $|a_1-a_2|>|a_1|$. Then for $0<\ka<|a_1|\we|a_2| $, we have that
\beg{align*}
\ka^2- |2a_1-a_2|\ka +a_1(a_1-a_2) & =(\ka-|a_1|)(\ka-|a_1-a_2|)\\
&\geq (\ka-|a_1|\we |a_2|)(\ka-|a_1-a_2|).
\end{align*}
Thus, by \eqref{in-sika}, 
$$p(\ka)-2\al\ka^2-\|\si\|_\infty^2> 0.$$
This together with \eqref{Ar0} yields that $r_0< \ka$. Taking into account \eqref{Ar0}, we get \eqref{sset3}.  Hence, \eqref{eq-ex1} has a stationary probability measure $\nu_1\in \sP^{M,1}_{a_1,\ka}$. It is clear that we can  replace $a_1$ by $a_2$. Thus \eqref{eq-ex1} has a stationary probability measure $\nu_2\in \sP^{a_1,1}_{ M,\ka}$. Finally, since
$$\ka<\ff {|a_1|\we|a_2|} 2\leq \ff 1 2 \left(\ff {|a_1-a_2|} 2\right),$$
the two probability measures $\nu_1$ and $\nu_2$ are distinct by Theorem \ref{thmN}.   

If $a_1=-a_2$ and $\si$ is a constant, it has been proved that \eqref{eq-ex1} has a symmetric stationary probability measure (see e.g. \cite{HT10a}):
\beg{align*}
\mu_0(\d x)=c_0\exp\left\{-\ff 2 {\si^{2}}\left(\be\left(\ff {x^4} 4-\ff {a_1^2x^2} 2\right)-\ff {\al x^2} 2\right) \right\}\d x,
\end{align*}
where $c_0$ is the normalization  constant so that $\mu_0$ is a probability measure.  Assuming $a_1>0$, it follows from \eqref{nunu} and $\ka< {a_1} /2$ that 
$$\nu_1([0,+\infty))=\nu_1(|\cdot-a_1| \geq a_1)\leq \ff {\nu_1(|\cdot-a_1|)} {a_1}\leq \ff {\ka} {a_1}<\ff 1 2,$$
and, similarly, $\nu_2((-\infty,0])< 1/2$.  Thus, $\nu_1$ and $\nu_2$ are not symmetric measures. Hence, \eqref{eq-ex1} has three stationary distributions.

\qed

\bigskip

\noindent{\bf Proof of Example \ref{ex-1}~}

According to \eqref{eq-ex2}, we have that
\beg{align*}
b(x,\mu)&=-\ff {\be} 2\left\{(x-a_1)|x-a_2|^2+(x-a_2)|x-a_1|^2\right\}\\
&\qquad -\int_{\R^d}\left(\al_1 |x-y|^2(x-y)+\al_2 (x-y)\right)\mu(\d y).
\end{align*}
It is a  routine to check that all the conditions of Corollary \ref{cor0}, and we omit the details. The rest is prove  \eqref{bgg}, \eqref{sset1} and \eqref{sset2} for 
\beg{align*}
g_1(r,r_1,r_2)&= 2\be\Big\{\left(\ff {|a_1-a_2|^2} 2+\ff {\al_1} {\be}r_1^2+\ff {\al_2}{\be}\right)r^2 -\left(\ff {\al_1} {\be}r_2^3+\ff {\al_2} {\be} r_1\right)r\Big\}-\|\si\|_\infty^2,\\
g_{2}(r,r_1) & =2\be\left\{ \left(1+\ff {\al_1} {\be}\right)r^{\ff 4 3}-3\left(|a_1-a_2|+\ff {\al_1} {\be} r_1\right)r \right\},\\
g (r,\mu)&=g_1(r,\mu)+g_2(r^3,\mu),~r\geq 0.
\end{align*}

It is clear that
\beg{align*}
&-\ff {\be} 2\left(|x|^2|x+a_1-a_2|^2+\<x+a_1-a_2,x\>|x|^2\right)\\
&\qquad = -\be\left(|x|^4+\ff 3 2\<a_1-a_2,x\>|x|^2+\ff {|a_1-a_2|^2} 2|x|^2\right)\\
&\qquad \leq -\be\left(|x|^4-\ff 3 2|a_1-a_2| |x|^3+\ff {|a_1-a_2|^2} 2|x|^2\right),
\end{align*}
and 
\beg{align*}
&-\al_1\int_{\R^d}|x+a_1-y|^2\<x+a_1-y,x\>\mu(\d y) -\al_2\int_{\R^d}\<x+a_1-y,x\>\mu(\d y)\\
&\qquad = -\al_1\int_{\R^d}\left(|x|^2+2\<x,a_1-y\>+|a_1-y|^2\right)\left(|x|^2+\<a_1-y,x\>\right)\mu(\d y)\\
&\qquad\quad  -\al_2\int_{\R^d}\left(|x|^2+\<a_1-y,x\>\right)\mu(\d y)\\
&\qquad \leq -\al_1\left(|x|^4+3|x|^2\int_{\R^d}\<x,a_1-y\>\mu(\d y)+|x|^2\int_{\R^d}|a_1-y|^2\mu(\d y)\right)\\
&\qquad\quad -\al_1\left(2\int_{\R^d}\<a_1-y,x\>^2\mu(\d y)+ \int_{\R^d}|a_1-y|^2\<a_1-y,x\>\mu(\d y)\right)\\
&\qquad\quad  -\al_2 |x|^2-\al_2\int_{\R^d}\<a_1-y,x\> \mu(\d y)\\
&\qquad \leq -\al_1\left(|x|^4-3 \mu(|\cdot-a_1| ) |x|^3+\mu(|\cdot-a_1|^2)|x|^2-\mu(|\cdot-a_1|^3)|x|\right)\\
&\qquad\quad -\al_2|x|^2+\al_2\mu(|\cdot-a_1|)|x|.
\end{align*}
Thus
\beg{align*}
2\<b(x+a_1,\mu),x\>&=- \be \left(|x|^2|x+a_1-a_2|^2+\<x+a_1-a_2,x\>|x|^2\right)\\
&\quad -2\al_1\int_{\R^d}|x+a_1-y|^2\<x+a_1-y,x\>\mu(\d y)\\
&\quad -2\al_2\int_{\R^d}\<x+a_1-y,x\>\mu(\d y)\\
& \leq  -2\be\Big\{\left(1+\ff {\al_1} {\be}\right)|x|^4-3\left(|a_1-a_2|+\ff {\al_1} {\be} \mu(|\cdot-a_1| )\right) |x|^3\\
&\quad +\left(\ff {|a_1-a_2|^2} 2+\ff {\al_1} {\be}(\mu(|\cdot-a_1|))^2+\ff {\al_2}{\be}\right)|x|^2\\
&\quad -\left(\ff {\al_1} {\be}\mu(|\cdot-a_1|^3)+\ff {\al_2} {\be} \mu(|\cdot-a_1| ) \right)|x|\Big\}.
\end{align*}
Then 
\beg{align*}
2\<b(x+a_1,\mu),x\>+\|\si(x+a_1,\mu)\|_{HS}^2\leq - g(|x|,\|\mu_{a_1}\|_1,\|\mu_{a_1}\|_3),~x\in\R^d.
\end{align*}
For $\mu$ with $\|\mu_{a_1}\|_1 <\ka_1$ and $\|\mu_{a_1}\|_3<\ka_2$. Since $\ff {\al_1} {\be}>1$, we have that
\beg{align*}
\th_1&=\ff {(1+\ff {\al_1} {\be})^{\ff 3 2}\left(\ff {\be} {\al_1}\right)^2} {2(3(4+\th_0))^{\ff 3 2}}\leq \ff {(  \ff {2\al_1} {\be})^{\ff 3 2}\left(\ff {\be} {\al_1}\right)^2} {2(3(4+\th_0))^{\ff 3 2}} < \ff 1 2\sq{\ff {\be} {\al_1}}.
\end{align*}
Let $\th_2=\ff {\ka_2} {|a_1-a_2|}$. Then 
$$\th_1 <\left(\ff {\be \th_1} {4\al_1} \right)^{\ff 1 3}=\th_2.$$ As a consequence, we have that $\ka_1<\ka_2$. It is clear that 
\beg{align*}
(\pp_r^2g) (r,r_1,r_2)&=2\be\Big\{12\left(1+\ff {\al_1} {\be}\right)r^2-18\left(|a_1-a_2|+\ff {\al_1} {\be} r_1\right) r\\
&\quad +2\left(\ff {|a_1-a_2|^2} 2+\ff {\al_1} {\be} r_1^2+\ff {\al_2}{\be}\right)\Big\}.
\end{align*}
By  \eqref{aba1a2}, we have that 
\beg{align*}
\th_1 < \left(\ff 2 {\sq{27}} \sq{(1+\ff {\al_1}{\be}) (1+2\th_0)}-1\right)\ff {\be} {\al_1}.
\end{align*} 
Then
\beg{align*}
\left(1+\ff {\al_1\th_1} {\be}\right)^2<\ff 4 {27} \left(1+\ff {\al_1} {\be} \right)(1+2\th_0).
\end{align*}
This implies that
\beg{align*}
\left(18\left(|a_1-a_2|+\ff {\al_1} {\be} \mu(|\cdot-a_1| )\right)\right)^2<48\left(1+\ff {\al_1} {\be}\right) \left( {|a_1-a_2|^2} +\ff {2\al_2}{\be}\right).
\end{align*}
Consequently,  $g(\cdot,r_1,r_2)$ is convex.  Thus there is $r_0>0$ such that 
\beg{align*}
&\left\{r\geq 0~|~g (r,r_1,r_2)\leq 0\right\}\\
&\qquad\subset \Big\{r\geq 0~\big|~ (1+\ff {\al_1} {\be})r^4-3( |a_1-a_2|+\ff {\al_1} {\be}\ka_1)r^3\\
&\qquad \qquad+\left(\ff {|a_1-a_2|^2} 2 +\ff {\al_2}{\be}\right)r^2-\left(\ff {\al_1} {\be}\ka_2^3+\ff {\al_2} {\be} \ka_1\right)r \leq \ff {\|\si\|_\infty^2} {2\be}\Big\}\\
&\qquad =  [0,r_0].
\end{align*}
Since \eqref{aba1a2}, we have that
Thus 
\beg{align*}
\ff {\th_1} 2-3\th_1^2+(1-\ff {2\al_1} {\be})\th_1^3>\ff {9\th_1} {20}>\ff {\th_1} 4.
\end{align*}
This together with $\ff 1 4\th_1= \ff {\al_1} {\be}\th_2^3$ and $\ka_1=\th_1|a_1-a_2|$ implies that 
\beg{align*}
\ff {|a_1-a_2|^2} 2\ka_1-3|a_1-a_2|\ka_1^2+(1-\ff {2\al_1} {\be})\ka_1^3> \ff {\al_1} {\be}\ka_2^3. 
\end{align*}
Moreover, it follows from \eqref{si-ka12} that 
\beg{align*}
&(1+\ff {\al_1} {\be})\ka_1^4-3(|a_1-a_2|+\ff {\al_1} {\be}\ka_1)\ka_1^3+\left(\ff {|a_1-a_2|^2} 2 +\ff {\al_2}{\be}\right)\ka_1^2-\left(\ff {\al_1} {\be}\ka_2^3+\ff {\al_2} {\be} \ka_1\right)\ka_1\\
&\qquad \geq  (1- \ff {2\al_1} {\be})\ka_1^4- 3|a_1-a_2|\ka_1^3 + \ff {|a_1-a_2|^2} 2  \ka_1^2- \ff {\al_1} {\be}\ka_2^3  \ka_1\\
&\qquad =\left( (1- \ff {2\al_1} {\be})\th_1^3- 3 \th_1^2 + \ff {\th_1} 4 \right)\th_1 |a_1-a_2|^4\\
&\qquad > \ff {\th_1^2} 5|a_1-a_2|^4>\ff {\|\si\|_\infty^2} {2\be}. 
\end{align*}
Hence, $r_0\leq \ka_1$. Hence, \eqref{sset1} holds.\\
For  $(r_1,r_2)\in\sD_{\ka_1,\ka_2}$, we have that 
\beg{align*}
& g_{2}(\ka_2^3,r_1) +\inf_{\ka_1\geq r\geq 0}g_1(r,r_1,r_2) \\
&\quad\geq  2\be\left\{ \left(1+\ff {\al_1} {\be}\right)\ka_2^{4}-3\left(|a_1-a_2|+\ff {\al_1} {\be}\ka_1\right)\ka_2^3 \right\}- \ff {\be\left(\ff {\al_1}{\be}\ka_2^3+\ff {\al_2} {\be}\ka_1\right)^2} {  |a_1-a_2|^2+ \ff {2\al_2} {\be} }-\|\si\|_\infty^2\\
&\quad=2\be\Big\{ \left(1+\ff {\al_1} {\be}\right)\th_2^{4}-3\left(1+\ff {\al_1} {\be}\th_1\right)\th_2^3 - \ff { \left(\ff {\al_1}{\be}\th_2^3+\th_0\th_1\right)^2} { 2+ 4\th_0 }\Big\}|a_1-a_2|^4-\|\si\|_\infty^2 .
\end{align*}
Taking into account that $ \ff {\be\th_1} {4\al_1}=  \th_2^3$, we have by \eqref{si-ka12} that
\beg{align*}
&  \left(1+\ff {\al_1} {\be}\right)\th_2^{4}-3\left(1+\ff {\al_1} {\be}\th_1\right)\th_2^3 - \ff { \left(\ff {\al_1}{\be}\th_2^3+\th_0\th_1\right)^2} { 2+ 4\th_0 } \\
&\quad \geq  \left\{\left(1+\ff {\al_1} {\be}\right)\left(\ff {\be } {4\al_1}\right)^{\ff 4 3}\th_1^{\ff 1 3}-\ff {3\be} {4\al_1}-\left(\ff 3 4+\ff {(\ff 1 4+\th_0)^2} {2+4\th_0}\right)\th_1\right\}\th_1\\
&\quad > \left\{\left(1+\ff {\al_1} {\be}\right)\left(\ff {\be } {4\al_1}\right)^{\ff 4 3}\th_1^{\ff 1 3}-\ff {3\be} {4\al_1}-\left(1+\ff {\th_0} 4 \right)\th_1\right\}\th_1\\
&\quad =  \ff {(1+\ff {\al_1} {\be})^{\ff 3 2}\left(\ff {\be}{\al_1}\right)^2} {24\sq3(1+\ff {\th_0} 4)^{\ff 1 2}}-\ff {3\be} {4\al_1}=\ff {(4+\th_0)\th_1} {2}-\ff {3\al_1}{4\be}\\
&\quad >\ff {\|\si\|_\infty^2} {2\be|a_1-a_2|^4}. 
\end{align*}
Thus  
$$ g_{2}(\ka_2^3,r_1)+g_1(r,r_1,r_2)>0,~(r_1,r_2)\in\sD_{\ka_1,\ka_2}.$$
Note that $g_2(\cdot,r_1)$ decreases first and then increase. Hence,
\beg{align*}
\bigcup_{(r_1,r_2)\in\sD_{\ka_1,\ka_2}}\left\{0 \leq r~\Big|~g_{2}(r^3,r_1)+\inf_{0\leq r \leq \ka_1}g_1(r ,r_1,r_2)\leq 0\right\}\subset [0,\ka_2].
\end{align*}
Therefore, \eqref{sset2} holds, and the proof is completed. 

\qed

\section{Existence results for singular case}

In this section, we investigate the existence of stationary distributions of distribution depended SDEs with singular drifts by using the Zvonkin transformation initialed in \cite{Zv}. Well-posedness results for \eqref{equ00} have been established by \cite{HW} recently. While only stationary distributions are consider in this paper, we use more weaker conditions on the coefficients, see (H4) and (H5) below. 

We denote $L^p$ the usual $L^p$-space on $\R^d$  and by $\|\cdot\|_p$ the $L^p$-norm. For $(\th,p)\in [0,2]\times (1,+\infty)$, we define $H^{\th,p}=(\1-\De)^{-\ff {\th} 2}(L^p)$ to be the usual Bessel potential space with norm
$$\|f\|_{\th,p}=\|(\1-\De)^{\ff {\th} 2}f\|_p,$$
where $\De$ is the Laplace operator on $\R^d$ and $\1$ is the identity operator.  Let $\chi\in C_c^\infty(\R^d)$ with $\1_{[|x|\leq 1]}\leq \chi\leq \1_{[|x|\leq 2]}$. We define  
$$\chi_r(x)=\chi(\ff x r),\qquad \chi_r^z(x)=\chi\left(\ff {x-z} r\right),~r>0,x,z\in\R^d.$$
We denote by $\tld H^{\th,p}$ the localized $H^{\th,p}$-space introduced in \cite{XXZZ}:
$$\tld H^{\th,p}:=\left\{f\in H^{\th,p}_{loc}(\R^d)~|~\|f\|_{\tld H^{\th,p}}:=\sup_z\|\chi_r^z f\|_{\th,p}<\infty\right\}.$$
In particular, we denote $\tld L^p=\tld H^{0,p}$. Given a probability $\mu$. We consider \eqref{equ0} in the following form 
\beg{align}\label{equ-Ph}
\d X_t^\mu=b_0(X_t^\mu,\mu)\d t+b_1(X_t^\mu,\mu)\d t+\si(X_t^\mu,\mu)\d W_t.
\end{align}
The drift term $b_0$ is regular, and $b_1$ is singular satisfying the following hypotheses. 

\beg{description}[align=left, noitemsep]
\item[(H4)] There exists $p>d$  such that  
$$\ka_0:=\sup_{\mu\in\sP^{1+r_2}}\|b_1(\cdot,\mu)\|_{\tld L^p}<\infty.$$
For every $n\geq 1$ and $M\geq 1$
\beg{align}\label{cb1Lp}
\lim_{ \nu\xrightarrow{w} \mu \text{~in~} \sP^{1+r_2}_{M}}\|(b_{1}(\cdot,\mu)-b_{1}(\cdot,\nu))\1_{[|\cdot|\leq n]}\|_{L^p}=0.
\end{align}
\end{description}
Instead of the locally Lipschitz as \eqref{b0-K}, we assume that $\si$ satisfies
\beg{description}[align=left, noitemsep]
\item[(H5)] For every $\mu\in\sP^{1+r_2}$, $\si(\cdot,\mu)$ is uniformly continuous on $\R^d$ and $\nn\si(\cdot,\mu)\in \tld L^p$.
\end{description}
Then we have the following theorem under some uniformly non-degenerate condition on the diffusion term $\si$.

\beg{thm}\label{thm2}
Assume that $b_0$ satisfies   (H1)-(H3) with setting $\si\equiv 0$ there and that \eqref{ineb-g} is replaced by the following stronger condition:
\beg{align}\label{ineb-gs}
|b_0(x,\mu)|\leq C_8(1+|x|^{r_1})+C_9\|\mu\|_{1+r_2}^{\ff {r_3r_1} {1+r_1}},~x\in\R^d,
\end{align} 
where $C_8$ is a positive constant independent of $\mu$;  $b_1$ satisfies (H4); $\si$ satisfies (H3), (H5) and  there are positive constants $\la_1,\la_2$ such that 
\beg{align}\label{node}
\la_1\leq (\si\si^*)(x,\mu)\leq \la_2,~x\in\R^d,\mu\in\sP^{1+r_2}.
\end{align}
If $r_3 \leq   1+r_1   $, and $C_1> C_3$ when $r_3 = 1+r_1 $, then \eqref{equ-Ph} has a stationary distribution.
\end{thm}

We finally give an example on the non-uniqueness of stationary distributions  for  distribution depended SDEs with a measurable bounded  drift.

\beg{exa}\label{ex-3}
Let $d=1$.  Consider the SDE  in Example \ref{ex-0} perturbed by a bounded drift: 
\beg{align}\label{eq-ex3}
\d X_t&=-\be(X_t-a_1)X_t(X_t-a_2)\d t-\al \int_{\R}(X_t-y)\sL_{X_t}(\d y)\d t\nonumber\\
&\quad +h(X_t,\sL_{X_t})\d t+\si(X_t,\sL_{X_t})\d W_t,
\end{align}
where the constants $\al,\be,a_1,a_2$  satisfy all the conditions of  Example \ref{ex-0}, $\si$ satisfies (H5) and \eqref{node}, and $h$ is a bounded measurable function on $\R\times\sP(\R)$ satisfying \eqref{cb1Lp}. If there is $\ka\in (0, (|a_1|\we|a_2|)/2)$ such that
\beg{align}\label{in-sika-n}
\|\si\|_\infty^2+\|h\|_\infty \ka< 2\be\ka ^2(\ka-|a_1-a_2|)(\ka-|a_1|\we |a_2|) ,
\end{align} 
then there exist two distinct stationary distributions  $\nu_1,\nu_2\in\sP^{1+r_2}$ such that
\beg{align*} 
\nu_1(|\cdot-a_1|)\leq \ka,\qquad \nu_2(|\cdot-a_2|)\leq \ka.
\end{align*}

\end{exa}

\subsection{Proofs of Theorem \ref{thm2} and Example \ref{ex-3}}

\noindent {\bf{Proof of Theorem \ref{thm2}.~~}} 

For given $\mu\in\sP^{1+r_2}(\R^d)$, under the assumptions of Theorem \ref{thm2}, it has been proved that \eqref{equ-Ph} has a unique strong solution and the associated semigroup $P_t^\mu$ has a unique invariant probability measure, which we also denote by $\cT_\mu$.  We still denote by $\cT$ the mapping on $\sP^{1+r_2}$.

Let $a(x,\mu)=(\si\si^*)(x,\mu)$, and let  $u_\mu$ be the solution of the following equation
\beg{align}\label{bPDE}
\ff 1 2 \trac\left(a(\cdot,\mu)\nn^2 u_\mu \right)(x)+(\nn_{b_1(x,\mu)} u_\mu)(x)=\la u_\mu(x)-b_1(x,\mu).
\end{align}
According to \cite[Theorem 2.1]{YZ} or \cite[Theorem 7.5]{XZ}, $u_\mu\in \tld H^{2,p}$ and 
\beg{align}\label{nnu}
\lim_{\la\ra+\infty}\sup_{\mu\in\sP^{1+r_2}}\left(\|u_\mu\|_\infty+\|\nn u_\mu\|_\infty\right)=0.
\end{align}
 Let $U_\mu(x)=x+u_\mu(x)$.  Then by the It\^o formula,  we get that
\beg{align}\label{UX}
 \d U_\mu(X_t^{\mu})  &=\left(\nn U_\mu b_{0 }(\cdot,\mu)+\la u_\mu\right) (X_t^\mu) \d t+(\nn U_\mu \si(\cdot,\mu))(X_t^\mu)\d W_t. 
\end{align}
By \eqref{nnu}, we choose $\la_0>0$ such that  for $\la>\la_0$, $\sup_{\mu\in\sP^{1+r_2}}\|\nn u_\mu\|_\infty<\ff 1 2$.  Then
\beg{align}\label{Uup}
\sup_{\mu\in\sP^{1+r_2}}\left(\|\nn (U_\mu)^{-1}\|_\infty\vee \|\nn U_\mu\|_\infty\right)&\leq \sup_{\mu\in\sP^{1+r_2}}\ff 1 {1-\|\nn u_\mu\|_\infty}\leq 2.
\end{align}
Thus $U^\mu$ is a diffemorphism on $\R^d$. 

\beg{lem}
Under the assumptions of Theorem \ref{thm2}, there exists $\la_0>0$ and for each $\la>\la_0$, there is $M_0>0$ such that for any $M\geq M_0$, $\sP^{1+r_2}_{M}$ is invariant under the mapping $\cT$.
\end{lem}
\beg{proof}
Let $Y_t^\mu=U_\mu(X_t^\mu)$. Then by \eqref{UX}, 
\beg{align}\label{eqY} 
 \d Y_t^\mu  &=\left(\nn U_\mu b_{0 }(\cdot,\mu)+\la u_\mu\right) (U_\mu^{-1}(Y_t^\mu)) \d t+(\nn U_\mu \si(\cdot,\mu))(U_\mu^{-1}(Y_t^\mu))\d W_t. 
\end{align}
It follows from \eqref{node}, \eqref{Uup} and (H5) that $(\nn U_\mu \si(\cdot,\mu))(U_\mu^{-1}(y))$ is continuous in $y$ and non-degenerate. Define a mapping $\cT\circ U^{-1}$ on $\sP^{1+r_2}$ as follows
$$(\cT\circ U^{-1})_\mu=\cT_\mu\circ U_\mu^{-1},~\mu\in\sP^{1+r_2}.$$  
Then $(\cT \circ U^{-1})_\mu$ is the invariant probability measure of \eqref{eqY}.

Next, we check  the coefficients of \eqref{eqY} is subject to (H1) and (H2) except \eqref{b0-K}. Then all the assumptions  of Lemma \ref{lem1}, Lemma \ref{lem2} and \eqref{lem3} hold. 

We first verify (H2) except \eqref{b0-K}. By \eqref{node} and \eqref{Uup}, it is clear that the diffusion term  $\nn U_\mu \si(\cdot,\mu)$ of \eqref{eqY} satisfies \eqref{ine-b-si00}:
\beg{align}\label{supUsi}\sup_{\mu\in\sP^{1+r_2},y\in\R^d}\|(\nn U_\mu \si(y,\mu))(U_\mu^{-1}(y))\|_{HS}^2\leq 4\la_2^2=:\tld C_6.
\end{align}
For every $x\in\R^d$
\beg{align}\label{ine-U0}
\left||x|-\|u\|_\infty\right|\leq |(U_\mu)^{-1}(x)|&=|x-u((U_\mu)^{-1}(x))|\leq  |x|+\|u_\mu\|_\infty. 
\end{align}
Then  the Jensen inequality yields that
\beg{align}\label{U-1low}
|U_\mu^{-1}(x)|^{1+r_1}&\geq (1-\ep_0)^{-r_1} |x|^{1+r_1}-\ep_0^{-r_1}\|u\|_\infty^{1+r_1},~\ep_0\in (0,1),\\
|U_\mu^{-1}(x)|^{ r_1}&\leq 2^{(r_1-1)^+}|x|^{r_1}+2^{(r_1-1)^+}\|u_\mu\|_\infty^{r_1}.\label{U-1up}
\end{align}
Inequality \eqref{U-1up} together with \eqref{ineb-gs} yields that 
\beg{align}\label{b0-si}
|b_0(U_\mu^{-1}(x),\mu)|\leq C_8\left[1+2^{(r_1-1)^+}(|x|^{r_1}+\|u_\mu\|^{r_1})\right]+C_9\|\mu\|_{1+r_2}^{\ff {r_3r_1} {1+r_1}}.
\end{align}
Hence, \eqref{ineb-g} holds.

We then verify (H1) for \eqref{eqY}. Since \eqref{supUsi}, we  focus on the drift term. By \eqref{ine-b-si0}, we have that 
\beg{align}\label{ineUb}
&\<(\nn U_\mu b_0(\cdot,\mu))(U_\mu^{-1}(x)),x\>+\la\<u_\mu(U_\mu^{-1}(x)),x\>\nonumber\\
&\qquad\leq \<b_0(U_\mu^{-1}(x),\mu),x\>+\|\nn u_\mu\|_\infty|b_0(U_\mu^{-1}(x),\mu)||x|+\la \|u_\mu\|_\infty|x|\nonumber\\
&\qquad\leq \<b_0(U_\mu^{-1}(x),\mu), U_\mu^{-1}(x)\>+\<b_0(U_\mu^{-1}(x),\mu), u_\mu(U_\mu^{-1}(x))\>\nonumber\\
&\qquad\quad  +\|\nn u_\mu\|_\infty|b_0(U_\mu^{-1}(x),\mu)||x|+\la \|u_\mu\|_\infty|x|\nonumber\\
&\qquad\leq -\ff {C_1}2|U_\mu^{-1}(x)|^{1+r_1}+\ff {C_2} 2+\ff {C_3} 2\|\mu\|_{1+r_2}^{r_3}\nonumber\\
&\qquad\quad+ \|u_\mu\|_\infty|b_0(U_\mu^{-1}(x),\mu)|+\|\nn u_\mu\|_\infty|b_0(U_\mu^{-1}(x),\mu)||x|+\la \|u_\mu\|_\infty|x|.
\end{align}
This and the H\"older inequality imply that for any positive $\ep_1,\ep_2,\ep_3$,
\beg{align}\label{bbu}
&\|u_\mu\|_\infty|b_0(U_\mu^{-1}(x),\mu)|+\|\nn u_\mu\|_\infty|b_0(U_\mu^{-1}(x),\mu)||x|\nonumber\\
&\quad \leq \left(\ff {r_1} {1+r_1}\ep_1^{\ff {1+r_1} {r_1}}+\|\nn u_\mu\|_\infty C_8 2^{(r_1-1)^+}+(1+r_1)^{-1}\ep_2^{1+r_1}+\ff {C_9\|\nn u_\mu\|_\infty} {1+r_2}\right)|x|^{1+r_1}\nonumber\\
&\qquad +C_8\|u_\mu\|_\infty^{1+r_1}2^{(r_1-1)^+}+(1+r_1)^{-1}\left(\ff {\|u_\mu\|_\infty C_8 2^{(r_1-)^+}} {\ep_1}\right)^{1+r_1}\nonumber\\
&\qquad +\ff {r_1} {1+r_1}\left(\ff {\|\nn u_\mu\|_\infty C_8(1+2^{(r_1-1)^+}\|u_\mu\|^{r_1})} {\ep_2}\right)^{\ff {1+r_1} {r_1}}+\ff {\|u_\mu\|_\infty C_9} {(1+r_1)\ep_2^{1+r_1}} \nonumber\\
&\qquad +\left(\ff {r_1 C_9\|u_\mu\|_\infty \ep_3^{\ff {1+r_1} {r_1} }} {1+r_1} +\ff {C_9\|\nn u_\mu\|_\infty r_1} {1+r_1}\right)\|\mu\|_{1+r_2}^{r_3}\nonumber\\
&\quad =: \bar C_1(\ep_1,\ep_2,\mu,\la)|x|^{1+r_1}+\bar C_2(\ep_1,\ep_2,\ep_3,\mu,\la)+\bar C_3(\ep_3,\mu,\la)\|\mu\|_{1+r_2}^{r_3}.
\end{align}
Putting this and \eqref{U-1low} into \eqref{ineUb}, and by using the H\"older inequality, we have that 
\beg{align*}
&\<(\nn U_\mu b_0(\cdot,\mu))(U_\mu^{-1}(x)),x\>+\la\<u_\mu(U_\mu^{-1}(x)),x\>\\
&\qquad  \leq -\left(\ff {C_1} {2(1-\ep_0)^{r_1}}-\bar C_1( \ep_1,\ep_2,\mu,\la)-\ff {\ep_4^{r_1+1}} {1+r_1} \right)|x|^{1+r_1}\\
&\quad\qquad +\bar C_2(\ep_1,\ep_2,\ep_3,\mu,\la)+\ff {\|u_\mu\|_\infty^{1+r_1}} {2\ep_0^{r_1}}+\ff {r_1(\la\|u_\mu\|_\infty)^{\ff {r_1} {1+r_1}}} {(1+r_1)\ep_4^{\ff {r_1} {1+r_1} }}\\
&\quad\qquad +\left(\ff {C_3} 2+\bar C_3(\ep_3,\mu,\la)\right)\|\mu\|_{1+r_2}^{r_3},
\end{align*}
where $\ep_4$ is any positive constant.  Hence, (H1) holds. Moreover, for given $\la>\la_0$ and $\ep_0, \ep_1,\ep_3, \ep_4$
\beg{align*}
\sup_{\mu\in\sP^{1+r_2}}\left\{(\bar C_1+\bar C_2+\bar C_3)(\ep_0,\ep_1,\ep_2,\ep_3,\mu,\la)+\ff {\|u_\mu\|_\infty^{1+r_1}} {2\ep_0^{r_1}}+\ff {r_1(\la\|u_\mu\|_\infty)^{\ff {r_1} {1+r_1}}} {(1+r_1)\ep_4^{\ff {r_1} {1+r_1} }}\right\}<\infty.
\end{align*}

Due to \eqref{nnu}, it is clear that 
$$\lim_{\la\ra+\infty}\sup_{\mu\in\sP^{1+r_2}}(\bar C_1( \ep_1,\ep_2,\mu,\la)+\bar C_3(\ep_3,\mu,\la))=0$$
As a consequence, in the case that $C_1>C_3$, we can choose $\ep_0, \ep_1,\ep_3, \ep_4$ small and a larger $\la_0>0$ such that for $\la>\la_0$ 
\beg{align}\label{C13}
\sup_{\mu\in\sP^{1+r_2}}\ff { \ff {C_3} 2+\bar C_3(\ep_3,\mu,\la)} {\ff {C_1} {2(1-\ep_0)^{r_1}}-\bar C_1(\ep_1,\ep_2,\mu,\la)-\ff {\ep_4^{r_1+1}} {1+r_1} } <\ff {C_3+1} {C_1+1}<1.
\end{align}
Up to now,   we have prove  that there exists $\la_0>0$ and for any $\la>\la_0$, the coefficients of \eqref{eqY} satisfy (H1)  and \eqref{inesi} with some  other constants  $\tld C_i,~i=1,\cdots,6$  independent of $\mu$.  Moreover,  $\tld C_4=\tld C_5=0$ by \eqref{supUsi}, and due to \eqref{C13}, we also have  that  
\beg{align}\label{C3/1}
\ff {\tld C_3} {\tld C_1}<\ff {C_3+1} {C_1+1}<1,~\text{when}~r_3= 1+r_1 .
\end{align} 
This implies that \eqref{ine-C-r} and \eqref{ine-C-r'} hold.

Though we can find an invariant subset of $\cT\circ U^{-1}$ in $\sP^{1+r_2}$ by  Lemma \ref{lem2},  
we need to  prove that $\sP^{1+r_2}_M$ with large enough $M$ is  an invariant subset of $\cT$ instead of $\cT\circ U^{-1}$. Following the proof of Lemma \ref{lem2}, we have by \eqref{ine-F-mu1} for $\cT\circ U^{-1}$ that
\beg{align*}
\|(\cT\circ U^{-1})_\mu\|_{1+r_2}^{ 1+r_1 }&\leq \ff {\tld C_2+\tld C_6(r_2-r_1)} {\tld C_1 } +\ff {\tld C_3 } {\tld C_1 } \|\mu\|_{1+r_2}^{r_3}.
\end{align*}
In the case that $r_3= {1+r_1} $, by \eqref{C3/1}, there exists $\tld M_0>0$ such that
\beg{align*}
\ff {\tld C_2+\tld C_6(r_2-r_1)} {\tld C_1 } +\ff {\tld C_3 } {\tld C_1 }M^{r_3}\leq \ff {C_3+1} {C_1+1} M^{ {1+r_1} },~M\geq  \tld M_0. 
\end{align*} 
In the case that $r_3< {1+r_1} $, it is clear that there exists $\tld{\tld{M_0}}>0$ such that  
\beg{align*}
\ff {\tld C_2+\tld C_6(r_2-r_1)} {\tld C_1 } +\ff {\tld C_3 } {\tld C_1 }M^{r_3}\leq \left(\ff {C_3+1} {C_1+1} \we \ff 1 2\right)M^{ {1+r_1}  },~M\geq \tld{\tld{M_0}}. 
\end{align*} 
 Consequently, there exists $c_0\in (0,1)$ such that for $M\geq  \tld M_0\vee \tld{\tld{M_0}}$, 
\beg{align*}
\| (\cT\circ U^{-1})_\mu\|_{1+r_2}^{ 1+r_1 }&\leq \ff {\tld C_2+\tld C_6(r_2-r_1)} {\tld C_1 } +\ff {\tld C_3 } {\tld C_1 }M^{r_3}\\
&\leq c_0M^{ {1+r_1} },~\mu\in\sP^{1+r_2}_{M}.
\end{align*}
This together with \eqref{ine-U0} implies that
\beg{align*}
\|\cT_\mu\|_{1+r_2}&=\left(\cT_\mu\circ U_\mu^{-1}(|U_\mu^{-1}(\cdot)|^{1+r_2})\right)^{\ff 1 {1+r_2}} \\
&\leq \left(\cT_\mu\circ U_\mu^{-1}(|\cdot|^{1+r_2})\right)^{\ff 1 {1+r_2}}+\|u_\mu\|_\infty^{\ff 1 {1+r_2}}\\
&\leq c_0^{\ff 1 {1+r_1}}M + \sup_{\mu\in\sP^{1+r_2}}\|u_\mu\|_\infty^{\ff 1 {1+r_2}},~\mu\in\sP^{1+r_2}_{M}.
\end{align*}
Thus, for 
$$M\geq M_0:= \tld M_0\vee\tld{\tld{M_0}}\vee \ff {\sup_{\mu\in\sP^{1+r_2}}\|u_\mu\|_\infty } { (1-c_0^{\ff 1 {1+r_1}})^{ (1+r_2)}},$$
we have that 
\beg{align*}
\|\cT_\mu\|_{1+r_2} \leq c_0^{\ff 1 {1+r_1}}M + (1-c_0^{\ff 1 {1+r_1}}) M = M^{\ff 1 {1+r_2}},~\mu\in\sP^{1+r_2}_{M}.
\end{align*}
Therefore, $\sP^{1+r_2}_{M}$ is invariant under $\cT$.

\end{proof}

We next prove that $\cT$ is weak sequence continuous  in $\sP^{1+r_2}_{M}$. Under the assumptions of Theorem \ref{thm2}, by applying Lemma \ref{lem3} to \eqref{eqY}, it is clear that $P_t^\mu$ is   $V$-uniformly exponential ergodic with some locally bounded $V$.  According to the proof of Lemma \ref{lem4}, the weak sequence continuity of $\cT$ follows from the lemma below. Let 
$$\ta_n=\inf\left\{t>0~|~ |X_t^{\nu,x}|+|X_t^{\mu,x}|>n \right\}.$$
\beg{lem}
Under the assumptions of Theorem \ref{thm2}, and let $\mu,\mu_m\in\sP^{1+r_2}_{M}$ such that $\mu_m\xrightarrow{w} \mu$. Then
\beg{align}\label{in-lim1}
\lim_{m\ra+\infty}\int_{|x|\leq \ff n 2}\E\left(\left|X_{t \we\ta_n}^{\mu,x}-X_{t \we\ta_n}^{\nu_m,x}\right|\right)\cT_{\nu_m}(\d x)=0,~t>0.
\end{align} 
\end{lem}
\beg{proof} Let $X_t^{n,\nu}$ be the solution of the following SDE
\beg{align}\label{equ-Phn}
\d X_t^{n,\nu}=(b_0\et_n)(X_t^{n,\nu},\nu)\d t+b_1(X_t^{n,\nu},\nu)\d t+\si(X_t^{n,\nu},\nu)\d W_t,
\end{align}
where $\et_n\in C^1_b(\R^d)$ is a cutoff function  with $\1_{[|x|\leq n]}\leq \et_n \leq \1_{[|x|\leq n+1]} $. By the pathwise uniqueness, 
\beg{align*}
X_t^{n,\mu}=X_t^{\mu},~~X_t^{n,\nu}=X_t^{\nu},~t<\ta_n
\end{align*}
Denote $b_{0n}=b_0\et_n$. Then
\beg{align*}
\d u^\mu(X_t^{n,\nu}) & =\left(\ff 1 2 \trac\left(a(\cdot,\nu)\nn^2 u^\mu \right)(X_t^{n,\nu})+(\nn_{b_1(\cdot,\nu)} u^\mu+\nn_{b_{0n}(\cdot,\nu)}u^\mu)(X_t^{n,\nu})\right)\d t\\
&\qquad  +(\nn_{\si(\cdot,\nu)\d W_t} u^\mu)(X_t^{n,\nu}).
\end{align*}
Then by using  \eqref{bPDE}, we have that 
 \beg{align*}
 \d U^\mu(X_t^{n,\nu})&=\d X_t^{n,\nu}+\d u^\mu(X_t^{n,\nu})\\
 & =\left( b_{0n}+b_1\right)(X_t^{n,\nu},\nu) \d t+\si(X_t^{n,\nu},\nu)\d W_t+\nn_{\si(X_t^{n,\nu},\nu)\d W_t} u^\mu(X_t^{n,\nu})\\
 &\qquad + \left(\ff 1 2 \trac\left(a(\cdot,\nu)\nn^2 u^\mu \right)\right)(X_t^{n,\nu})+\left(\nn_{b_1(\cdot,\nu)} u^\mu  +\nn_{b_{0n}(\cdot,\nu)} u^\mu\right)(X_t^{n,\nu})\d t\\
&=\left( b_{0n}+b_1\right)(X_t^{n,\nu},\nu) \d t+\si(X_t^{n,\nu},\nu)\d W_t+\nn_{\si(X_t^{n,\nu},\nu)\d W_t} u^\mu(X_t^\nu)\\
&\qquad +  \left(\ff 1 2 \trac\left(a(\cdot,\mu)\nn^2 u^\mu \right)(X_t^{n,\nu})+(\nn_{b_1(X_t^{n,\nu},\mu)} u^\mu)(X_t^{n,\nu})\right)\d t\\
&\qquad + \ff 1 2 \trac\left(\left(a(\cdot,\nu)-a(\cdot,\mu)\right)\nn^2 u^\mu \right)(X_t^{n,\nu})\d t+(\nn_{b_1(\cdot,\nu)-b_1(\cdot,\mu)} u^\mu)(X_t^{n,\nu})\d t\\
&\qquad+(\nn_{b_{0n}(\cdot,\nu)}u^\mu)(X_t^{n,\nu}) \d t  \\
&=\left( b_{0n}(\cdot,\nu)+\nn_{b_{0n}(\cdot,\nu)}u^\mu +b_1(\cdot,\nu)-b_1(\cdot,\mu)+\la u^\mu\right)(X_t^{n,\nu}) \d t\\
&\qquad + \ff 1 2 \trac\left(\left(a(\cdot,\nu)-a(\cdot,\mu)\right)\nn^2 u^\mu \right)(X_t^{n,\nu})\d t+(\nn_{b_1(\cdot,\nu)-b_1(\cdot,\mu)} u^\mu)(X_t^{n,\nu})\d t\\
&\qquad  +\si(X_t^{n,\nu},\nu)\d W_t+\nn_{\si(X_t^{n,\nu},\nu)\d W_t} u^\mu(X_t^{n,\nu})\\
&= \left\{\nn U^\mu b_{0n}(\cdot,\nu)+\nn U^\mu \left(b_1(\cdot,\nu)-b_1(\cdot,\mu)\right)+\la u^\mu \right\}(X_t^{n,\nu})\d t\\
&\qquad + \ff 1 2 \trac\left(\left(a(\cdot,\nu)-a(\cdot,\mu)\right)\nn^2 u^\mu \right)(X_t^{n,\nu})\d t  +\nn U^\mu(X_t^{n,\nu})\si(X_t^{n,\nu},\nu)\d W_t.
 \end{align*} 
In particular, 
\beg{align*}
 \d U^\mu(X_t^{n,\mu})  &=\left(\nn U^\mu b_{0n}(\cdot,\mu)+\la u^\mu\right) (X_t^{n,\mu}) \d t+\nn U^\mu (X_t^{n,\mu})\si(X_t^{n,\mu},\mu)\d W_t. 
\end{align*} 
Then
\beg{align*}
\d \left(U^\mu(X_t^{n,\mu})-U^\mu(X_t^{n,\nu})\right)&=\left\{\nn U^\mu(X_t^\mu)b_{0n}(X_t^\mu,\mu)-\nn U^\mu(X_t^\nu)b_{0n}(X_t^\nu,\nu) \right\}\d t\\
&\qquad + \nn U^\mu(X_t^{n,\nu})\left(b_1(X_t^{n,\nu},\mu)-b_1(X_t^{n,\nu},\nu)\right)\d t\\
&\qquad +\la \left(u^\mu(X_t^{n,\mu})-u^\mu(X_t^{n,\nu}) \right)  \d t\\
&\qquad + \ff 1 2 \trac\left(\left(a(\cdot,\mu)-a(\cdot,\nu)\right)\nn^2 u^\mu \right)(X_t^{n,\nu})\d t \\
&\qquad +\left((\nn U^\mu\si(\cdot,\mu)) (X_t^{n,\mu})-(\nn U^\mu\si(\cdot,\nu))(X_t^{n,\nu})\right)\d W_t.
\end{align*}
Then
\beg{align*}
&\d |U^\mu(X_t^{n,\mu})-U^\mu(X_t^{n,\nu})|^2\\
&\quad =2\left\<(\nn U^\mu b_{0n}(\cdot,\mu))(X_t^{n,\mu})-(\nn U^\mu b_{0n}(\cdot,\nu))(X_t^{n,\nu}), U^\mu(X_t^{n,\mu})-U^\mu(X_t^{n,\nu})\right\> \d t\\
&\quad\quad + 2\left\<\left(\nn U^\mu \left(b_1(\cdot,\mu)-b_1(\cdot,\nu)\right)\right)(X_t^{n,\nu}), U^\mu(X_t^{n,\mu})-U^\mu(X_t^{n,\nu})\right\>\d t\\
&\quad\quad +2\la \left\< u^\mu(X_t^{n,\mu})-u^\mu(X_t^{n,\nu}), U^\mu(X_t^{n,\mu})-U^\mu(X_t^{n,\nu})\right\> \d t\\
&\quad\quad +   \left\<\trac\left(\left(a(\cdot,\mu)-a(\cdot,\nu)\right)\nn^2 u^\mu \right)(X_t^{n,\nu}), U^\mu(X_t^{n,\mu})-U^\mu(X_t^{n,\nu})\right\>\d t \\
&\quad\quad +\|\nn U^\mu (X_t^{n,\mu})\si(X_t^{n,\mu},\mu)-\nn U^\mu(X_t^{n,\nu})\si(X_t^{n,\nu},\nu)\|_{HS}^2\d t\\
&\quad\quad +2\left< U^\mu(X_t^{n,\mu})-U^\mu(X_t^{n,\nu}),\left((\nn U^\mu \si(\cdot,\mu))(X_t^{n,\mu})-(\nn U^\mu\si(\cdot,\nu))(X_t^{n,\nu})\right)\d W_t\right\>\\
&\quad =(I_1+I_2+I_3+I_4+I_5)\d t+\d M_t,
\end{align*}
where
$$M_t   =2 \int_0^t \left< U^\mu(X_s^{n,\mu})-U^\mu(X_s^{n,\nu}),\left((\nn U^\mu \si(\cdot,\mu))(X_s^{n,\mu})-(\nn U^\mu\si(\cdot,\nu))(X_s^{n,\nu})\right)\d W_t\right\>.$$
We choose $\la>0$ large enough such that 
$$\ff 1 2\leq\|\nn U^\mu\|_\infty\leq \ff 3 2.$$
By the Krylov's estimate, see e.g.  \cite[Theorem 5.6]{XZ} or \cite[Theorem 3.1]{YZ}, the distributions of $X_t^{n,\mu}$ and $X_t^{n,\nu}$ are absolutely w.r.t. the Lebesgue measure for almost every $t>0$. Combining this with  \eqref{b0-K}, $\nn U^\mu\in C_b(\R^d;\R^d\otimes\R^d)$,  we have that
\beg{align*}
I_1 &=2\left\<\left(\nn U^\mu(X_t^{n,\mu})-\nn U^\mu(X_t^{n,\nu})\right)b_{0n}(X_t^{n,\mu},\mu), U^\mu(X_t^{n,\mu})- U^\mu(X^{n,\nu})\right\>\\
&\qquad +2\left\<\nn U^\mu(X_t^{n,\nu})\left(b_{0n}(X^{n,\mu},\mu)-b_{0n}(X^{n,\nu},\mu)\right), U^\mu(X_t^{n,\mu})- U^\mu(X_t^{n,\nu})\right\>\\
&\qquad +2\left\<\nn U^\mu(X_t^{n,\nu})\left(b_{0n}(X^{n,\nu},\mu)-b_{0n}(X^{n,\nu},\nu)\right), U^\mu(X_t^{n,\mu})- U^\mu(X_t^{n,\nu})\right\>\\
&\leq \ff {3C} 2 \left(\sup_{|x|\leq n} |b_0(x)|\right) \left(|\cM_1(\nn^2 u^\mu)(X_t^{n,\mu})|+|\cM_1(\nn^2 u^\mu)(X_t^{n,\nu})| \right) |X_t^{n,\mu}-X_t^{n,\nu}|^2 \\
&\qquad + \ff {9} 2 K_n |X_t^{n,\mu}-X_t^{n,\nu}|^2+\ff 9 2|X_t^{n,\mu}-X_t^{n,\nu}|\sup_{|x|\leq n}|b_0(x,\mu)-b_0(x,\mu)|\\
&\leq \tld C_n\left(|\cM_1(\nn^2 u^\mu )(X_t^{n,\mu})|+|\cM_1(\nn^2 u^\mu )(X_t^{n,\nu})|+1 \right)|X_t^{n,\mu}-X_t^{n,\nu}|^2\\
&\qquad +\sup_{|x|\leq n}|b_0(x,\mu)-b_0(x,\mu)|^2,~t\leq \ta_n,
\end{align*} 
holds for some $\tld C_n>0$. \\
For $I_2$, 
\beg{align*}
I_2 & \leq \ff {9} 2\left|b_1(X_t^{n,\nu},\mu)-b_1(X_t^{n,\nu},\nu)\right|\cdot |X_t^{n,\mu}-X_t^{n,\nu}|\\
&\leq \ff {9} 2 \left|b_1(X_t^{n,\nu},\mu)-b_1(X_t^{n,\nu},\nu)\right|^2+\ff 9 2|X_t^{n,\mu}-X_t^{n,\nu}|^2\\
&=\ff {9} 2 \left|b_1(X_t^{n,\nu},\mu)-b_1(X_t^{n,\nu},\nu)\right|^2+\ff 9 2|X_t^{n,\mu}-X_t^{n,\nu}|^2.
\end{align*}
For $I_3$,
\beg{align*}
I_3\leq    \ff {3\la } 2 |X_t^{n,\mu}-X_t^{n,\nu}|^2.
\end{align*}
For $I_4$, 
\beg{align*}
I_4& \leq \ff 3 2\|a(X_t^{n,\nu},\mu)-a(X_t^{n,\nu},\nu) \|_{HS}\|\nn^2 u^\mu(X_t^{n,\nu})\|_{HS}|X_t^{n,\mu}-X_t^{n,\nu}|\\
&\leq  3\|\si\|_\infty\|\nn^2 u^\mu(X_t^{n,\nu})\|_{HS}\|\si(X_t^\nu,\mu)-\si(X_t^\nu,\nu)\|_{HS} |X_t^{n,\mu}-X_t^{n,\nu}|\\
&\leq \ff 3 2\|\nn^2 u^\mu(X_t^{n,\nu})\|_{HS}^2 |X_t^{n,\mu}-X_t^{n,\nu}|^2+\ff 3 2\|\si\|_\infty^2\|\si(X_t^{n,\nu},\mu)-\si(X_t^{n,\nu},\nu)\|_{HS}^2.
\end{align*}
For $I_5$, 
\beg{align*}
I_5&\leq  \|\si\|_\infty^2\|\nn U^\mu (X_t^{n,\mu})-\nn U^\mu(X_t^{n,\nu})\|_{HS}^2+\ff 9 4\|\si(X_t^{n,\mu},\mu)- \si(X_t^{n,\nu},\nu)\|_{HS}^2\\
&\leq C\|\si\|_\infty^2\left(\|\cM_1(\nn^2 u^\mu ) (X_t^\mu)\|_{HS}+\|\cM_1(\nn u^\mu )(X_t^\nu)\|_{HS}^2\right)|X_t^\mu-X_t^\nu|^2\\
&\qquad +C\left(\|\cM_1(\nn\si(\cdot,\mu))(X_t^{n,\mu})\|_{HS}^2+\|\cM_1(\nn\si(\cdot,\mu))(X_t^{n,\nu})\|_{HS}^2\right)|X_t^{n,\mu}-X_t^{n,\nu}|^2\\
&\qquad +\ff 9 2\|\si(X_t^{n,\nu},\mu)-\si(X_t^{n,\nu},\nu)\|_{HS}^2.
\end{align*}
Hence,  there exist  $C_1'>0$,
\beg{align*}
  |U^\mu(X_{t\we\ta_n}^\mu)-U^\mu(X_{t\we\ta_n}^\nu)|^2 & =|U^\mu(X_{t\we\ta_n}^{n,\mu})-U^\mu(X_{t\we\ta_n}^{n,\nu})|^2\\
&\leq \int_0^{t\we\ta_n}|X_s^{n,\mu}-X_s^{n,\nu}|^2\d A_{n,s}^{\mu,\nu} +t\sup_{|x|\leq n}|b_0(x,\mu)-b_0(x,\mu)|^2\\
 &\qquad + \ff 9 2\int_0^{t\we\ta_n}|b_1(X_s^{n,\nu},\mu)-b_1(X_s^{n,\nu},\nu)|^2\d s\\
 &\qquad + C_1't\sup_{|x|\leq n}\|\si(x,\mu)-\si(x,\nu)\|_{HS}^2+ M_{t\we\ta_n}
\end{align*}
where
\beg{align*}
&A_{n,t}^{\mu,\nu}  = \int_0^{t }\Big\{\tld C_n\left(1+ \|\cM_1(\nn^2 u^\mu ) (X_s^{n,\mu})\|_{HS}^2+\|\cM_1(\nn^2 u^\mu )(X_s^{n,\nu})\|_{HS}^2\right)\\
& +C_2'\left( \|\cM_1(\nn\si(\cdot,\mu))(X_s^{n,\mu})\|_{HS}^2+(\|\nn^2u^\mu \|_{HS}^2+\|\cM_1(\nn\si(\cdot,\mu))\|_{HS}^2)(X_s^{n,\nu})\right) \Big\}\d s
\end{align*} 
for some positive constants $\tld C_n$, which depends on $n$,  and $C_2' $. By stochastic Gronwall's lemma, we have that
\beg{align*}
&\left(\E | X_{t\we\ta_n}^\mu -  X_{t\we\ta_n}^\nu | \right)^2 \leq 4 \left(\E |U^\mu(X_{t\we\ta_n}^\mu)-U^\mu(X_{t\we\ta_n}^\nu)|\right)^2\\
&\qquad\leq 2c_pC_2'\left(\E e^{\ff p {p-1} A_{n,t\we\ta_n}^{\mu,\nu}}\right)^{\ff {1-p} p}t\Big\{ \sup_{|x|\leq n}\|\si(x,\mu)-\si(x,\nu)\|_{HS}^2\\
&\qquad\qquad +\sup_{|x|\leq n}|b_0(x,\mu)-b_0(x,\mu)|^2+\E\int_0^{t\we\ta_n}|b_1(X_s^{n,\nu},\mu)-b_1(X_s^{n,\nu},\nu)|^2\d s\Big\}.
\end{align*}
It is clear that
\beg{align*}
\int_0^{t\we\ta_n}|b_1(X_s^{n,\nu},\mu)-b_1(X_s^{n,\nu},\nu)|^2\d s&=\int_0^{t\we\ta_n}|b_1(X_s^{n,\nu},\mu)-b_1(X_s^{n,\nu},\nu)|^2\1_{[|X_s^{n,\nu}|\leq n]}\d s\\
&=\int_0^{t\we\ta_n}|b_{1n}(X_s^{n,\nu},\mu)-b_{1,n}(X_s^{n,\nu},\nu)|^2\d s
\end{align*}
where $b_{1n}(x,\nu)=b_1(x,\nu)\1_{[|x|\leq n]}$. For $q>\ff {2p} {p-d}$, it follows from the Krylov estimate that
\beg{align*}
\E\int_0^{t\we\ta_n}|b_{1n}(X_s^{n,\nu},\mu)-b_{1n}(X_s^{n,\nu},\nu)|^2\d s&\leq C t^{\ff 2 q }\|b_{1n}(\cdot,\mu)-b_{1n}(\cdot,\nu)\|_{\tld L^p}^2\\
&\leq C t^{\ff 2 q }\|b_{1n}(\cdot,\mu)-b_{1n}(\cdot,\nu)\|_{ L^p}^2
\end{align*}
where the constant $C$ depends on $p,q,t,n,\|b_{1n}(\cdot,\nu)\|_{L^p}$. For fixed $\mu\in \sP^{1+r_2}_{M}$ and a sequence $\nu_m\in\sP^{1+r_2}_{M}$ such that $\nu_m\xrightarrow{w}\mu$,  we have that
$$\lim_{m\ra+\infty}\|b_{1n}(\cdot,\nu_m)-b_{1n}(\cdot,\mu)\|_{L^p}=0$$
Then the sequence $\{\|b_{1n}(\cdot,\nu_m)\|_{L^p}\}_{m\geq 1}$ is bounded for any $n\geq 1$. 
Moreover,  this together with  Khasminskii's estimate (see e.g. \cite[Lemma 3.5]{XZ}) for $X_t^{n,\nu_m}$ and $X_t^{n,\mu}$ implies that
$$\sup_{m}\E e^{\ff p {p-1} A_{n,t\we\ta_n }^{\mu,\nu_m}}\leq \sup_{m}\E e^{\ff p {p-1} A_{n,t }^{\mu,\nu_m}}<\infty.$$
Hence
\beg{align*}
&\int_{|x|\leq \ff n 2}\E\left(\left|X_{t \we\ta_n}^{\mu,x}-X_{t \we\ta_n}^{\nu_m,x}\right|\right)\cT_{\nu_m}(\d x)\\
&\quad \leq C_n t\Big\{ \sup_{|x|\leq n}\left(\|\si(\cdot,\mu)-\si(\cdot,\nu_m)\|_{HS}+ |b_0(\cdot,\mu)-b_0(\cdot,\nu_m)|\right)(x)\\
&\qquad\quad + t^{\ff 1 q}\|b_{1n}(\cdot,\mu)-b_{1n}(\cdot,\nu_m)\|_{L^p}\Big\}.
\end{align*}
Therefore \eqref{in-lim1} holds. 

\end{proof}

So far, we have proved that there exists $M_0>0$ such that for any $M\geq M_0$, $\cT$ maps $\sP^{1+r_2}_{M}$ into $\sP^{1+r_2}_{M}$ and $\cT$ is  continuous in $\sP^{1+r_2}_{M}$ w.r.t. to the Kantorovich-Rubinstein distance. Therefore,  the assertion of this theorem follows from the Schauder fixed point theorem as proving in Theorem \ref{thm1}. 

\qed

\bigskip

\noindent{\bf Proof of Example \ref{ex-3}}

With $g(r,\mu)$ in the proof of Example \ref{ex-0} replaced by $g(r,r_1)-\|h\|_\infty r$, it follows from the proof of Example \ref{ex-0} that $\bigcup_{\mu\in\sP^{a_1,1}_{ M,\ka}} A_\mu^1\subset [0,r_0]$ for some $r_0\geq 0$. Let $p(r)$ be the polynomial defined in \eqref{pp}. Then it follows from \eqref{in-sika-n} that
\beg{align*}
& p(\ka)-2\al\ka^2-\|h\|_\infty \ka-\|\si\|_\infty^2\\
&\quad \geq 2\be\ka^2(\ka-|a_1|\we |a_2|)(\ka-|a_1-a_2|)-\|h\|_\infty \ka-\|\si\|_\infty^2\\
&\quad > 0.
\end{align*}
Hence $r_0\leq \ka$. Replacing $a_1$ by $a_2$ we have the same consequence. Therefore, there exist two distinct stationary distributions  $\nu_1\in \sP^{a_1,1}_{ M,\ka},\nu_2\in\sP^{a_2,1}_{ M,\ka}$.

\qed

\bigskip
 
\noindent\textbf{Acknowledgements}

\medskip

The  authors was supported by the National Natural Science Foundation of China (Grant No. 11901604, 11771326).



\begin{thebibliography}{99}



\bibitem{CarDe}
	R. Carmona, F. Delarue, Probabilistic Theory of Mean Field Games with Applications I, Springer, 2018.

\bibitem{CGPS}
	Carrillo, J.A., Gvalani, R.S., Pavliotis, G.A. and A. Schlichting, Long-Time Behaviour and Phase Transitions for the Mckean–Vlasov Equation on the Torus. Arch Rational Mech Anal 235, 635–690 (2020).  
	
\bibitem{Chen}
	Chen, M.-F., From Markov Chains to Non-Equilibrium Particle Systems, (2nd Ed.), World Scientific, 2004.
	
\bibitem{CMV}	
	J. A. Carrillo, R. J. McCann, C. Villani, Kinetic equilibration rates for granular media and related equations: entropy dissipation and mass transportation estimates, Rev. Mat. Iberoam. 19(2003), 971–1018.

\bibitem{ChPa}
	Chayes, L., Panferov, V., The McKean-Vlasov equation in finite volume. {\it J. Stat. Phys.} {\bf 138(1–3)}, 351--380, 2010

\bibitem{Daw}
	D. A. Dawson, Critical dynamics and fluctuations for a mean-field model of cooperative behavior. {\it J. Stat. Phys.} {\bf 31(1)}, 29--85, 1983
	
	
	
\bibitem{Deim}
	K. Deimling, Nonlinear functional analysis, Springer-Verlag Berlin Heidelberg, 1985. 	
	
	
\bibitem{DuTu}
	M. H. Duong, J. Tugaut, Stationary solutions of the Vlasov-Fokker-Planck equation: Existence, characterization and phase-transition, {\it Applied Mathematics Letters}  {\bf 52}, (2016), 38--45
	


\bibitem{EGZ}
	A. Eberle, A. Guillin, R. Zimmer, Quantitative Harris-type theorems for diffusions and McKean-Vlasov processes, Trans. Amer. Math. Soc. 371 (2019), 7135–7173.
	
\bibitem{FenZ}
	Feng, S. and Zheng, X. G., Solutions of a class of non-linear Master equations, {\it Stoch. Proc. Appl.} {\bf 43},(1992), 65-84. 






\bibitem{GoMa}
	B.~Goldys and B.~Maslowski, Exponential ergodicity for stochastic reaction-diffusion equations.
Stochastic partial differential equations and applications, XVII, 115--131, Lect. Notes Pure Appl.
Math., 245, Chapman Hall/CRC, Boca Raton, FL, 2006.






	
	


\bibitem{HT10a}
	S. Herrmann, J. Tugaut, Non-uniqueness of stationary measures for self-stabilizing processes, {\it stoch. Proc. Appl.} {\bf 120} (2010), 1215--1246.


\bibitem{HT10b}
	S. Herrmann and J. Tugaut, Stationary measures for self-stabilizing processes: Asymptotic analysis in the small noise limit, Electron. J. Probab. 15 (2010), pp. 2087–2116.
	
\bibitem{HW}
	X. Huang, F.-Y. Wang, Distribution dependent SDEs with singular coefficients, to appear at {\it Stoch. Proc. Appl.}
	
\bibitem{HRW}
	X. Huang, P. Ren, F.-Y. Wang, Distribution Dependent Stochastic Differential Equations, arXiv: 0212.13656.


















	
	





	
 


 

 

 



\bibitem{Mal}
	F. Malrieu,  Convergence to equilibrium for granular media equations and their Euler schemes, {\it Ann.  Appl. Probab.} {\bf 13(2)}, (2003)  540--560.



\bibitem{McK}
	H.~P.~McKean, Jr.,  Propagation of chaos for a class of non-linear parabolic equations. In: Stochastic Differential Equations (Lecture Series in Differential Equations, Session 7, Catholic Univ.), pp. 41--57. Air Force Office Sci. Res., Arlington, VA, 1967

\bibitem{Mel}
	S. M\'el\'eard, Asymptotic behaviour of some interacting particle systems; McKean-Vlasov and Boltzmann models, Probabilistic models for nonlinear partial differential equations (Montecatini Terme, 1995), Lecture Notes in Math., vol. 1627, Springer, Berlin, 1996, pp. 42--95. 
	
\bibitem{RenWa}
	P. Ren, F.-Y. Wang, Exponential convergence in entropy and Wasserstein distance for McKean-Vlasov SDEs, arXiv:2010.08950. 

\bibitem{Szn}
	A.-S., Sznitman, Topics in propagation of chaos. In: \'Ecole d'\'Et\'e de Probabilit\'es de Saint-Flour XIX-1989, volume 1464 of Lecture Notes in Math., pp. 165--251. Springer, Berlin, 1991


\bibitem{Tam}
	Y. Tamura, On asymptotic behaviors of the solution of a nonlinear diffusion equation. {\it J. Fac. Sci. Univ. Tokyo Sect. IA Math.} {\bf 31(1)}, (1984) 195--221. 
	
	
\bibitem{Tug10}	
	J. Tugaut, Convergence to the equilibria for self-stabilizing processes in double well landscape, Ann. Probab. (2010).  

\bibitem{Tug13}
	J. Tugaut, Self-stabilizing processes in multi-wells landscape in  $\R^d$-convergence, {\it Stoch. Proc. Appl.} {\bf 123} (2013) 1780--1801. 

\bibitem{Tug14a}
	J. Tugaut, Phase transitions of McKean–Vlasov processes in double-wells landscape, {\it Stochastics}, 86 (2) (2014), pp. 257-284

\bibitem{Tug14b}	
	J. Tugaut, Self-stabilizing processes in multi-wells landscape in $\R^d$-Invariant probabilities, {\it J. Theoret. Probab.}, 27 (1) (2014), pp. 57-79
	
	
	
\bibitem{Wan18}
	F.-Y. Wang, Distribution dependent SDEs for Landau type equations, Stoch. Proc. Appl. 128(2018), 595–621.
	
\bibitem{Wan21}
	F.-Y. Wang, Exponential Ergodicity for Fully Non-Dissipative McKean-Vlasov SDEs, arXiv: 2101.12562.

\bibitem{XXZZ}
	P. Xia, L. Xie, X. Zhang and G. Zhao, Lq(Lp)-theory of stochastic differential equations, arXiv: 1908.01255.
  



\bibitem{XZ}
	L. Xie, X. Zhang, Ergodicity of Stochastic Differential Equations with Jumps and Singular Coefficients, arXiv preprint arXiv:1705.07402.
	



\bibitem{YZ}
	C. Yuan, S.-Q. Zhang, A Zvonkin's transformation for stochastic differential equations with singular drift and related applications, arXiv: 1910.05903v5.
	




	








\bibitem{Zv}
	A.~K. Zvonkin, A transformation of the phase space of a diffusion process that will remove the drift, {\it Mat. Sb.}, {\bf 93} (135) (1974), 129--149.

\end{thebibliography}
\end{document}